\documentclass[12pt,leqno]{article}
\usepackage{amssymb}
\usepackage{hyperref}
\newcommand\qed{\hfill$\sqcap\kern-7.5pt\hbox{$\sqcup$}$}
\newcommand{\CC}{\mathbb{C}}

\newcommand{\ZZ}{\mathbb{Z}}
\newcommand{\RR}{\mathbb{R}}
\newcommand{\Ll}{\mathbb{L}}
\newcommand{\R}{\mathbb{R}}
\newcommand{\N}{\mathbb{N}}
\newcommand{\Sp}{\mathbb{S}}

\newcommand{\Ww}{\mathbb{W}}
\newtheorem{theo}{Theorem}
\newtheorem{prop}[theo]{Proposition}
\newtheorem{lem}[theo]{Lemma}
\newtheorem{cor}[theo]{Corollary}
\newtheorem{rem}[theo]{Remark}

\newcommand{\beqn}{\begin{equation}}
\newcommand{\eeqn}{\end{equation}}
\newcommand{\bear}{\begin{eqnarray}}
\newcommand{\eear}{\end{eqnarray}}
\newcommand{\bean}{\begin{eqnarray*}}
\newcommand{\eean}{\end{eqnarray*}}

\newcommand{\Sph}{\mathbb{S}}
\newcommand{\Aa}{\mathcal{A}}

\newcommand{\Cc}{\mathcal{C}}

\newcommand{\EE}{\mathcal{E}}
\newcommand{\FF}{\mathcal{F}}

\newcommand{\HH}{\mathcal{H}}

\newcommand{\LL}{\mathcal{L}}
\newcommand{\OO}{\mathcal{O}}

\newcommand{\Rr}{\mathcal{R}}

\newcommand{\Ss}{\mathcal{S}}

\newcommand{\chid}{{\chi_{_\delta}}}

\newcommand{\e}{{\varepsilon}}

\newcommand{\eps}{\varepsilon}

\begin{document}

\title{Stability, convergence to the steady state and elastic limit 
 for the Boltzmann equation for diffusively excited granular media} 

\date{}

\author{S. {\sc Mischler}$^1$, C. {\sc Mouhot}$^2$}

\footnotetext[1]{CEREMADE, Universit\'e Paris IX-Dauphine,
Place du Mar\'echal de Lattre de Tassigny, 75775 Paris, France. 
E-mail: \texttt{mischler@ceremade.dauphine.fr}}

\footnotetext[2]{CEREMADE, Universit\'e Paris IX-Dauphine,
Place du Mar\'echal de Lattre de Tassigny, 75775 Paris, France. 
E-mail: \texttt{mouhot@ceremade.dauphine.fr}}

\maketitle

\begin{abstract} 
We consider a space-homogeneous gas of {\it inelastic hard spheres}, with a {\it diffusive term} 
representing a random background forcing (in the framework of so-called 
{\em constant normal restitution coefficients} $\alpha \in [0,1]$ for the inelasticity). 
In the physical regime of a small inelasticity (that is $\alpha \in [\alpha_*,1)$ for some constructive 
$\alpha_* \in [0,1)$) we prove uniqueness of the stationary solution for 
given values of the restitution coefficient $\alpha \in [\alpha_*,1)$, the mass 
and the momentum, and we give various results on the linear stability and 
nonlinear stability of this stationary solution. 
\end{abstract}

\textbf{Mathematics Subject Classification (2000)}: 76P05 Rarefied gas
flows, Boltzmann equation [See also 82B40, 82C40, 82D05], 
76T25 Granular flows [See also 74C99, 74E20].

\textbf{Keywords}: Inelastic Boltzmann equation; granular gases; random forcing; 
hard spheres; stationary solution; uniqueness; stability; small inelasticity; 
elastic limit; degenerated perturbation; spectrum. 

\bigskip
\tableofcontents

\vspace{0.3cm}

\section{Introduction }
\setcounter{equation}{0}
\setcounter{theo}{0}


\subsection{The model and main result}

We consider the steady states of the spatially homogeneous inelastic Boltzmann equation for hard spheres with 
thermal bath forcing. More precisely, we consider a  gas which is described by the distribution density 
of particles $f= f(v) \ge 0$ with velocity $v \in \RR^N$ ($N \ge 2$) and such that $f$ satisfies the stationary equation
\bear  \label{eqBol1}
Q_\alpha(f,f) + \tau \, \Delta_v f  = 0 \quad\hbox{ in }\quad  \RR^N,
\eear
with a constant $\tau >0$, and given mass and vanishing momentum:
\bear
\label{eqBol2}
 \int_{\R^N} f \, dv = \rho \in (0,\infty), \quad \int_{\RR^N} f \, v \, dv = 0, 
\eear
\smallskip
The term $\mu \, \Delta_v \, f$, with constant $\mu \in (0,\infty)$, represents the effect of the heat bath. 
The quadratic collision operator $Q_\alpha(f,f)$ models the interaction of
particles by means of inelastic binary collisions with a constant normal restitution coefficient $\alpha \in [0,1)$ (see~\cite{GPV**,BGP**,MMRI,MMII}) which preserves  mass and momentum
but dissipates kinetic energy. We define the collision operator by its action on test functions, or 
{\it observables}. Taking $\psi = \psi(v)$ to be a suitably regular test function, 
we introduce the following weak formulation of the collision operator 
\beqn  \label{Qinel}
\int_{\R^N} 
Q_\alpha(g,f) \, \psi \, dv =  \int \!\! \int \!\!\int_{\R^N \times \R^N \times \Sp^{N-1}} 
b \, |u| \, g_* \, f \, (\psi' - \psi ) \, d\sigma \, dv \, dv_*, 
\eeqn
where we use the shorthand notations $f:= f(v)$,  $g_*:=g(v_*)$, $\psi':=\psi(v')$, etc.
Here and below $u=v-v_*$ denotes the relative velocity and $v',v'_*$ denotes the possible 
post-collisional velocities (which encapsule the inelasticity of the collision operator 
in terms of $\alpha$). They are defined by
\beqn\label{vprimvprim*}
\quad\quad   v' = {w \over 2} +  {u' \over 2}, \quad 
v'_*= {w \over 2} - {u' \over 2}, 
\eeqn
with
$$
w = v+v_*, \qquad u' = \left( {1- \alpha \over 2} \right) \, u + \left( {1+\alpha \over 2} \right) \, |u| \, \sigma.
$$

We also introduce the notation $\hat{x} = x/|x|$ for any $x \in \R^N$, $x\not= 0$.
The function $b=b(\hat u \cdot \sigma)$ in (\ref{Qinel}) is (up to a multiplicative factor)
the {\em differential collisional cross-section}. We assume that 
  \begin{equation}\label{hypb1}
  b \mbox{  is Lipschitz, non-decreasing and convex on } (-1,1) 
  \end{equation}
and that 
  \begin{equation}\label{hypb2}
  \exists \, b_m,b_M \in (0,\infty) \quad \mbox{ s.t. } \quad 
     \forall \, x \in [-1,1], \quad b_m \le b(x) \le b_M. 
   \end{equation}   
Note that the ``physical'' cross-section for hard spheres is given by
(see \cite{GPV**,Cerci?})
  \beqn\label{HScs} 
  b(x) = b'_0\, (1-x)^{-{N-3 \over 2}}, \quad b'_0 \in (0,\infty),
  \eeqn
so that it fulfills the above hypothesis (\ref{hypb1},\ref{hypb2}) when $N=3$. These hypothesis are 
needed in the proof of moments estimates (see~\cite[Proposition~3.2]{MMRI} 
and~\cite[Proposition~3.1]{MMII}).

We also define the symmetrized (or polar form of the) bilinear collisional operator $\tilde Q_\alpha$ by setting
 \beqn  \label{Qinelsym}
 \left\{
 \begin{array}{l}
 \displaystyle{
 \int_{\R^N} \tilde Q_\alpha(g,h) \, \psi \, dv = {1 \over 2} 
 \int \!\! \int \!\!\int_{\R^N \times \R^N \times \Sp^{N-1}} 
 b \, |u| \, g_* \, h \, \Delta_\psi \, d\sigma\,dv \, dv_*, }\vspace{0.3cm} \\
 \displaystyle{ \hbox{with} \quad \Delta_\psi = \left( \psi' + \psi'_* - \psi - \psi_*\right) }.
 \end{array}
 \right.
 \eeqn
In other words, $ \tilde Q_\alpha(g,h) =  ( Q_\alpha(g,h) +  Q_\alpha(h,g))/2$. 
The formula~(\ref{Qinel}) suggests the natural splitting 
$Q_\alpha = Q^+ _\alpha - Q^- _\alpha$ 
between gain and loss part. The loss part $Q^- _\alpha$ can be 
defined in strong form noticing that 
 $$
 \langle Q^-_\alpha(g,f) , \ \psi \rangle =  \int \!\! \int \!\!\int_{\R^N \times \R^N \times \Sp^{N-1}} 
 b \, |u| \, g_* \, f \, \psi  \, d\sigma \, dv \, dv_* =:  \langle f \, L(g) , \ \psi \rangle,
 $$
where $\langle \cdot, \cdot \rangle$ is the usual scalar product in $L^2$ and 
$L$ is the convolution operator
 \beqn\label{defL}\qquad
 L(g) (v) = (b_0 \,  |\cdot| * g )(v) = b_0 \, \int_{\R^N} g(v_*) \, |v-v_*| \, dv_*,
\,\,\,\hbox{with}\,\,\, b_0 = \int_{S^{N-1}} b(\sigma_1) \, d\sigma. 
 \eeqn
In particular note that $L$ and $Q^- _\alpha=Q^-$ are indeed independent of the 
normal restitution coefficient $\alpha$. 
\smallskip

As explained in~\cite{GPV**}, the operator~(\ref{Qinel})
preserves mass and momentum, and since the Laplacian also does so, 
the mass and momentum of a stationary solution can be prescribed. 
However energy is not preserved neither by the collisional operator (which tends to cool down the gas) 
nor by the diffusive operator (which harms it up). 
Competition between these two effects gives rise to a balance equation 
(obtained after having multiplied equation (\ref{eqBol1}) by $|v|^2$ and integrated) which reads as follows
\beqn \label{eqdiffEE}
(1-\alpha^2) \, D_\EE(f)  = \mu  \, 2 \, N \, \rho. 
\eeqn
The {\em energy dissipation functional} is given by
\[
D_\EE (f):= b_1 \int\!\!\int_{\RR^N \times \RR^N} f \, f_* \, |u|^3 \, dv \, dv_*, 
\]
where  $b_1$ is (up to a multiplicative factor) the angular momentum defined by 
\beqn \label{defb1}
b_1 :=  {1 \over 8} \int_{\Sp^{N-1}} \left( 1-(\hat{u}\cdot\sigma) \right) \, 
b(\hat{u}\cdot\sigma) \, d\sigma .
\eeqn
In order to establish (\ref{eqdiffEE}) we have used (\ref{Qinelsym}) and the elementary computation 
$$
\Delta_{|\cdot|^2} (v,v_*,\sigma) = - {1 - \alpha^2 \over 4} \, \, (1 -(\hat{u}\cdot\sigma)) \, |u|^2.
$$

\smallskip
Existence and qualitative properties of the steady solutions as well of the 
solutions to the associated evolution equation (Cauchy theory) was done in~\cite{GPV**} 
(see also \cite{BGP**}). Among others, it is proved the following in these papers:

\medskip
\begin{theo} \label{theo:exist} (\cite[Theorem 5.2, Lemma 7.2]{GPV**}, \cite[Theorem 1]{BGP**}) 
For any given inelastic coefficient $\alpha \in (0,1)$ and any given mass $\rho \in (0,\infty)$ 
there exists at least one solution $F \in \Ss(\R^N)$ to the stationary equation (\ref{eqBol1})-(\ref{eqBol2}) 
with  mass $\rho$ and vanishing momentum which furthermore satisfies
\beqn\label{upperlowerbdd}
F(v) \ge a_1 \, e^{-a_2 \, |v|^{3/2}}  \,\,\, \forall \, v \in \R^N
\quad\hbox{and}\quad \int_{\R^N} F(v) \, e^{a_3 \, |v|^{3/2}} \, dv < \infty,
\eeqn
for some constants $a_1, a_2, a_3 \in (0,\infty)$. .  
\end{theo}

Our main result gives a partial answer concerning uniqueness of solutions obtained by Theorem~\ref{theo:exist}:

\begin{theo} \label{theo:uniq}
There is some constructive $\alpha_* \in (0,1)$  such that for any $\alpha \in [\alpha_*,1]$, and any given mass 
$\rho \in (0,\infty)$ the stationary equation (\ref{eqBol1})-(\ref{eqBol2}) admits a unique solution $F$ 
with mass $\rho$ and vanishing momentum in the class of functions given by the existence Theorem~\ref{theo:exist}. 
\end{theo}


\subsection{Rescaled variables and elastic limit $\alpha \to 1$}

Let us start with a remark. For any solution $F$ to the equation~(\ref{eqBol1}) and any constant $\lambda \in (0,\infty)$, the rescaled function $g$ associated to $f$ by the relation
$$
g(v) = \lambda^N \, f (\lambda \,  v ).
$$
satisfies 
\beqn\label{eqresc}
Q_\alpha(g,g) + \lambda^{-3} \, \tau \, \Delta_v g = 0.
\eeqn
In order to prove (\ref{eqresc}) we have just use the homogeneity properties
$Q_\alpha(g,g)(v) = \lambda^{N-1} \, Q_\alpha(f,f)(\lambda v)$ and $\Delta_v \, g = \lambda^{N+2} \, (\Delta_v f) (\lambda \, v)$. 

\smallskip
This elementary remark shows that we may choose $\tau \in (0,\infty)$ arbitrarily  in equation~(\ref{eqBol1})  and we now make the choice 
\beqn\label{taualpha}
\tau = \tau_\alpha = \rho \, (1-\alpha), 
\eeqn
and denote by  $F_\alpha$ a solution to equation~(\ref{eqBol1})  with  $\tau = \tau_\alpha$:
\beqn\label{eqBolmua}
Q_\alpha(F_\alpha,F_\alpha) + \tau_\alpha \, \Delta F_\alpha = 0 \quad\hbox{in}\quad \R^N.
\eeqn
At a formal level, it is immediate that with this choice of scaling, 
in the elastic limit $\alpha \to1$, the equations (\ref{eqBol1})-(\ref{eqBol2}) become
  \begin{equation} \label{profileeq2}
  \left\{ 
  \begin{array}{l} \displaystyle 
  Q_1(F_1,F_1)  = 0  \,\,\hbox{ in }\,\,  \RR^N, \vspace{0.3cm}  \\ \displaystyle 
  \int_{\R^N} F_1 \, dv = \rho, \quad \int_{\RR^N} F_1 \, v \, dv = 0, 
     \quad 0 \le  F_1  \in \Ss(\R^N).  
  \end{array}
  \right. 
  \end{equation}
Moreover, from (\ref{eqdiffEE}), one gets
  \beqn\label{EnergyDEGe}
  2 \, N \,  (1-\alpha) \, \rho^2 - (1-\alpha^2) D_\EE(F_\alpha)  = 0.
  \eeqn
Dividing the above equation by $(1-\alpha)$ and passing to the limit $\alpha \to 1$, one obtains 
  \beqn\label{profilee3}
N \,   \rho^2 -  D_\EE(F_1)  = 0.
  \eeqn
It is straightforward that the only 
function satisfying the constraints (\ref{profileeq2}) and (\ref{profilee3}) 
is the Maxwellian function
  \beqn\label{Maxlim}
  \bar F_1 := M_{\bar\theta_1} = M_{\rho, 0, \bar\theta_1} 
  \eeqn
where, for any $\rho,\theta >0$, $u \in \R^N$, the function $M_{\rho,u,\theta}$ 
denotes the  Maxwellian with mass $\rho$, momentum $u$ and temperature 
$\theta$ given  by  
\beqn\label{defMrut}
M_{\rho,u,\theta} (v) :=  
  {\rho \over (2\pi \theta)^{N/2}} \, e^{- {|v-u|^2 \over 2 \theta}},
\eeqn
and where the temperature $\bar\theta_1 \in (0,\infty)$ is given by 
(we recall that $b_1$  is defined in~(\ref{defb1}))
  \beqn \label{tempSS}
  \bar\theta_1 =  {1 \over 2} \,  \frac{ N^{2/3} }{ b_1^{2/3} } \, 
         \left( \int_{\RR^N} M_{1,0,1}(v) \, |v|^3 \, dv \right)^{-2/3}.
  \eeqn
For instance in dimension $N=3$ we obtain $\bar\theta_1 = (3^2 \, \pi)^{1/3}/(2^{10} \, b_1^2)^{1/3}$.
Moreover,  in the particular case of the hard-spheres cross-section (\ref{HScs}) 
in dimension $3$, we find $b_1 = b'_0 (4 \pi)/3$ and therefore 
$\bar\theta_1 = 3^{4/3}/(2^{14}  \, \pi \, {b'_0}^2)^{1/3}$.


\subsection{Physical and mathematical motivation}



In this short paper we shall not review the physical and mathematical on the kinetic 
theory for granular gases. Let us only refer to some key references where the reader
can find all the desired details: for a detailed physical introduction we refer 
to~\cite{BPlivre,Cerci?}; and for a short mathematical introduction 
see~\cite{CVreviewgran}; finally apart from the many works on ``pseudo-maxwell molecules'' 
(see the bibliography of~\cite{CVreviewgran}) let us only mention the works 
directly related to inelastic hard spheres (freely cooling or in thermal bath): 
\cite{GPV**,BGP**,MMRI,MMII,MMGranular3} (in particular see the introductions of these 
papers for some mathematical discussions about the models and their connection 
to phyics). This paper is inspired from the 
models and questions arised in~\cite{GPV**} and the new mathematical tools 
developed in~\cite{MMGranular3}.


\subsection{Notation}\label{MMR:subsec:not}

Throughout the paper we shall use the notation
$\langle \cdot \rangle = \sqrt{1+|\cdot|^2}$.
We denote, for any $p \in [1,+\infty]$, $q \in \R$ and weight function $\omega:\R^N \to \R_+$, the 
weighted Lebesgue space $L^p_q(\omega)$ by 
     \[ 
     L^p_q(\omega) := \left\{f: \RR^N \mapsto \RR \hbox{ measurable }; \; \;
     \| f \|_{L^p_q(\omega)}   < + \infty \right\},
     \]
with, for  $p < +\infty$,
   \[ \| f \|_{L^p _q (\omega)} = \left[ \int_{\R^N} |f (v)|^p \, \langle v
        \rangle^{pq} \, \omega(v) \, dv \right]^{1/p} \]
and, for $p = +\infty$, 
    \[ \| f \|_{L^\infty _q (\R^N)} = \sup_{v \in \R^N} |f (v)| \, 
         \langle v \rangle^{q} \omega(v). \]
We shall in particular use the exponential weight functions 
\beqn\label{defdem}
m  = m_{s,a}(v) := e^{ -a \, |v|^s} \quad\hbox{for}\quad a  \in (0,\infty), \,\, s \in (0,1),
\eeqn    
or a smooth version $m(v):= e^{- \zeta(|v|^2)}$ with  $\zeta \in C^\infty$ is a positive 
function such that $\zeta(r) =  r^{s/2}$ for any $r \ge 1$,  with $s \in (0,1)$.

In the same way, the weighted Sobolev space $W^{k,p} _q (\omega)$ ($k \in \N$)
is defined by the norm
\[ 
\| f \|_{W^{k,p} _q (\omega)} =  \left[ \sum_{|s| \le k}
       \left\| \partial^s f (v) \right\|_{L^p _q (\omega)} ^p  \right]^{1/p}, 
\]
and as usual in the case $p=2$ we denote $H^k_q(\omega) =  W^{k,2} _q (\omega)$. The weight $\omega$ 
shall be omitted when it is $1$. 
Finally, for $f \in L^1_{2k}$, with $k \ge 0$, we introduce the following 
notation for the homogeneous moment of order $2k$
$$
{\bf m}_k(f) := \int_{\R^N} f \, |v|^{2 \, k} \, dv, 
$$
and we also denote by  $\rho(f) = {\bf m}_0(f)$ the mass of $f$, $\EE(f) = {\bf m}_1(f)$ 
the energy of $f$ and by $\theta(f) = \EE(f)/(\rho(f) \, N)$ the temperature associated 
to $f$ (when the distribution $f$ has $0$ mean). For any $\rho,\EE \in (0,\infty)$, 
$u \in \R^N$ we then introduce the subsets of $L^1$ of functions of given mass, 
mean velocity and energy
\bean
\Cc_{\rho,u}&:= &\left\{ h \in L^1_1; \, \int_{\R^N} h \, dv 
= \rho, \,\,  \int_{\R^N} h \, v \, dv = \rho \, u \right\}, \\
\Cc_{\rho,u,\EE} &:=& \left\{ h \in L^1_2; \, \int_{\R^N} h \, dv 
= \rho, \,\,  \int_{\R^N} h \, v \, dv 
= \rho \, u, \,\,  \int_{\R^N} h \, |v|^2 \, dv = \EE \right\}.
\eean
For any (smooth version of) exponential weight function $m$ and any $k,q \ge 0$ we introduce the Banach spaces
\beqn \label{defLL}
\Ll^1 (m^{-1}) = L^1(m^{-1}) \cap \Cc_{0,0} \quad\hbox{and}\quad
\Ww^{k,1}_q (m^{-1}) = W^{k,1}_q(m^{-1})  \cap \Cc_{0,0}.
\eeqn


\subsection{Summary of the results}

\smallskip
Our results, that we state now, deal with the evolution equation
\beqn\label{eqresca}
  \frac{\partial f}{\partial t} = Q_\alpha(f,f) + \tau_\alpha\, \Delta f, \quad 
  f(0,\cdot) = f_{\mbox{\scriptsize{in}}} \in \Cc_{\rho,0},
\eeqn
where $\tau_\alpha$  is defined in (\ref{taualpha}),  and with the associated stationary equation (\ref{eqBolmua}) with 
$F \in \Cc_{\rho,0}$.

\medskip
\begin{theo} \label{theo:uniq}
There is some constructive $\alpha_* \in (0,1)$  
such that for $\alpha \in [\alpha_*,1]$, and any given mass $\rho \in (0,\infty)$, we have:
\begin{itemize}
\item[(i)] For any $\tau>0$, the equation (\ref{eqBol1}) admits a unique non-negative 
stationary solution with  mass $\rho$ and vanishing momentum. 
We denote by $\bar F_\alpha$ the stationary solution obtained 
by fixing $\tau = \tau_{\alpha}$ (defined by (\ref{taualpha})). 
\item[(ii)] Let define $\bar F_1 = M_{\rho,0,\bar\theta_1}$ the Maxwellian  distribution with mass $\rho$, 
momentum $0$ and ``{\em diffusive thermodynamical temperature}'' $\bar\theta_1$ defined in~(\ref{tempSS}).
The path of stationary solutions $\alpha \to \bar F_\alpha$  parametrized by the normal restitution coefficient 
is $C^1$ from $[\alpha_*,1]$ into $W^{k,1} \cap L^1(e^{a \, |v|})$ for any $k \in \N$ and some $a \in (0,\infty)$.
\item[(iii)] For any $\alpha \in [\alpha_*,1]$, the linearized collision operator 
\beqn\label{defLa}
h \mapsto \LL_\alpha \, h := 2 \, \tilde Q_\alpha (\bar F_\alpha, h) + \tau_\alpha \, \Delta h
\eeqn
is well-defined and closed on $\Ll^1(m^{-1})$ for any exponential weight function $m$ with exponent $s\in(0,1)$ 
(defined in (\ref{defdem})).  
Its spectrum decomposes between a part 
which lies in the half-plane $\{ \mbox{{\em Re}} \, \xi \le \bar \mu \}$ for some constructive 
$\bar \mu < 0$, and some remaining discrete eigenvalue $\mu _\alpha$. This eigenvalue is real negative and satisfies 
\beqn\label{expanmua}
\mu _\alpha = - {3 \over \bar\theta_1} \, \rho \,  (1-\alpha) + \OO(1-\alpha)^2 \quad\hbox{when}\quad \alpha \to 1.
\eeqn
The associated eigenspace is of dimension $1$ and then denoting by $\phi_\alpha = \phi_\alpha(v)$ the unique associated 
eigenfunction such that $\| \phi_\alpha \|_{L^1_2} = 1$ and $\phi_\alpha(0) < 0$, there holds 
$\phi_\alpha \in \Ss(\R^N)$ (with bounds of regularity independent of $\alpha$) and 
\beqn\label{phiaTOphi1}
\phi_\alpha \rightarrow \phi_1 := c_0 \, \big( |v|^2 - N \, \bar\theta_1 \big) \, \bar F_1
\quad\hbox{ as }\quad  \alpha \to 1, 
\eeqn
where $c_0$ is the positive constant such that $\| \phi_1 \|_{L^1_2} = 1$. 
Finally one has constructive decay estimates on the semigroup 
associated to this spectral decomposition in this Banach space (see the 
key Theorem~\ref{thLalpha}Ê and the following point).

\item[(iv)] The stationary solution $\bar F_\alpha$ is globally attractive on bounded 
subsets of $L^1_3$ under some smallness condition on the inelasticity in the following sense. 
For any $\rho,\EE_0, M_0 \in (0,\infty)$ there exists $\alpha_{**} \in (\alpha_*,1)$,  
$C_* \in (0,\infty)$ and $\eta \in (0,1)$, such that for any initial datum  satisfying 
$$ 
0 \le f_{\mbox{\scriptsize{{\em in}}}}  \in L^1_3 \cap \Cc_{\rho,0,\EE_0}, 
\qquad \| f_{\mbox{\scriptsize{{\em in}}}} \|_{L^1_3} \le M_0,
$$
the solution $g$ to~(\ref{eqresca}) satisfies 
\beqn \label{cvgirreg} 
\| f_t - \bar F_\alpha \|_{L^1 _2}Ê\le e^{(1-\eta) \, \mu _\alpha \, t}.
\eeqn  

\item[(v)] Moreover, under smoothness condition on the initial datum one may 
prove a more precise asymptotic decomposition, and construct a Liapunov functional 
for the equation (\ref{eqresca}). 
More precisely, there exists $k_* \in \N$ and, for any exponential weight $m$ as defined in (\ref{defdem}) 
and any $\rho,\EE_0, M_0 \in (0,\infty)$, there exists $\alpha_{**} \in (\alpha_*,1)$ and 
a constructive functional $\HH : H^{k_*} \cap L^1(m^{-1}) \to \RR$ such that, first, 
for any initial datum 
$0 \le f_{\mbox{\scriptsize{{\em in}}}}  \in H^{k_*} \cap L^1(m^{-1}) \cap \Cc_{\rho,0,\EE_0}$ satisfying 
$$ 
\| f_{\mbox{\scriptsize{{\em in}}}} \|_{H^{k_*} \cap L^1(m^{-1})} \le M_0,
$$
the solution $f$ to~(\ref{eqresca}) satisfies 
\beqn\label{cvgcess}
f(t, \cdot) = \bar F_\alpha + c_\alpha (t) \, \phi_\alpha + r_\alpha(t, \cdot), 
\eeqn
with $c_\alpha (t) \in \R$ and $r_\alpha(t, \cdot) \in L^1_2(\R^N)$ such that 
\beqn\label{cvgcess2}
|c_\alpha (t)| \le  C_* \, e^{\mu _\alpha \, t}, \qquad
\| r_\alpha(t, \cdot ) \|_{L^1_2} \le C_* \, e^{(3/2) \, \mu_\alpha \, t}.
\eeqn 
And second when the initial datum satisfies additionally 
$$
f_{\mbox{\scriptsize{{\em in}}}} \ge M_0^{-1} \, e^{-M_0 \, |v|^8},
$$ 
the solution satisfies also
$$ 
t \mapsto \HH(g(t, \cdot)) \quad\hbox{is strictly decreasing }
$$
(up to reach the stationary state $\bar F_\alpha)$.
\end{itemize}
\end{theo}

\begin{rem} 

All the constants appearing in this theorem are contructive,  which means that they can be made explicit, 
and in particular that the proof does not use any compactness argument. Unless otherwise mentioned, 
these constants will depend on $b$, on the dimension $N$, and on some bounds on the initial datum 
but never on the inelasticity parameter $\alpha \in (0,1]$. 
All the other remarks made in \cite[Section~1.6]{MMGranular3} also apply here.
\end{rem}

\subsection{Method of proof and plan of the paper}

The general structure of the proof is inspired from~\cite{MMGranular3}. 
We have packed into the appendix the technical results from 
this previous paper which are used here also, as well as a technical 
result of lower bound on the diffusive inelastic Boltzmann equation. 
Therefore the overall structure of the proof is likely to be more 
visible than in~\cite{MMGranular3}. In the whole, we mainly prove here the few points which differs  from ~\cite{MMGranular3} (due to the replacement of the anti-drift term by a diffusive term) and refer to ~\cite{MMGranular3} for more
details.

The first main idea of our method is to consider the rescaled equations 
(\ref{eqBolmua})-(\ref{taualpha}) with an inelasticity dependent diffusion 
coefficient $\tau_\alpha$ which exactly ``compensates'' the loss of elasticity 
of the collision operator (in the sense that it compensates its loss of kinetic energy). 
This scaling allows to prove uniform bounds according to $\alpha$ for the family 
of stationary solutions $F_\alpha$ to the equation (\ref{eqBolmua}) (recall that in 
this scaling, the diffusion is evanescent in the elastic limit). 

The second main idea consists in decoupling the variations along the 
``energy direction'' and its ``orthogonal direction''. 
This decoupling makes it possible to identify the limit of different objects as 
$\alpha \to 1$ (among them the limit  of $F_\alpha$). 

The third main idea is to use systematically the knowledges on the elastic limit problem, 
once it has been identified thanks to the previous arguments. 
In particular we use the spectral study of the linearized problem and the dissipation 
entropy-entropy inequality for the elastic problem. 
This allows to argue by perturbative method. 
Let us emphasize that this perturbation is singular in the classical sense 
because of the addition of a (limit vanishing) 
second-order derivative operator, but also because of the gain of one more conservative quantity at the limit. 

\smallskip
In Section 2, we use the regularity properties of the collision operator in order to establish on the one hand 
that the family $(F_\alpha)$ is bounded in $H^\infty \cap L^1(m^{-1})$ uniformly according to the 
inelastic parameter $\alpha$ (the key argument being the use of the entropy functional which 
provides uniform lower bound on the energy of $F_\alpha$) and on the other hand that the difference 
of two stationary solutions in any strong norm may be bounded by the difference of these ones in weak norm 
(the key idea is a bootstrap argument). This last point shall allow to deal with the loss of derivatives and weights 
in the operator norms used in the sequel of the paper. 

\smallskip
In Section 3, we prove that $F_\alpha \to \bar F_1$  when $\alpha \to 1$ with explicit 
``H\"older'' rate. The cornerstone of the proof is again the decoupling of the variation 
$F_\alpha - \bar F_1$ between the ``energy direction'' and its ``orthogonal direction''.

\smallskip
Finally in Section 4, we prove uniqueness of the profile $\bar F_\alpha$ for small inelasticity 
by a variation around the implicit function theorem, in Section 5 we prove results 
on the localization of the spectrum of the linearized equation for small inelasticity, 
and in Section 6 we prove some (semi)-global stability results by combining the previous 
linearized study with entropy production estimates. 


\section{Estimates on the steady states} 
\setcounter{equation}{0}
\setcounter{theo}{0}


In this section we prove various regularity and decay estimates on the stationary solutions  
(or the differences of stationary solutions), uniform as $\alpha \to 1$, 
which shall be useful in the sequel. 

\subsection{Uniform estimates on the steady states}


For any $\alpha \in (0,1)$ we consider $\FF_\alpha$ the (not empty) set of all the solutions given by 
Theorem~\ref{theo:exist} of the inelastic diffusive stationary Boltzmann equation (\ref{eqBolmua}) 
with inelasticity coefficient $\alpha$, given mass $\rho \in (0,+\infty)$ and finite energy. 
More precisely, we define $\FF_\alpha$ as the following set of functions 
$$
\FF_\alpha := \Big\{ F \in  \Ss(\R^N) \quad \hbox{satisfying} \quad  
(\ref{eqBolmua}), \, (\ref{eqBol2}), \,(\ref{upperlowerbdd}\Big\}.
$$
For some fixed $\alpha_0\in (0,1)$, we also define 
$$
\FF  = \cup_{  \alpha \in [\alpha_0,1) }  \FF_\alpha.
$$ 

We show that for any stationary solution $F_\alpha \in \FF$ the decay estimates, 
the pointwise lower bound and the regularity estimates can be made uniform 
according to the inelasticity coefficient $\alpha \in [\alpha_0,1)$. Let us emphasize 
once again that the choice of the rescaling parameter $\tau_\alpha = \rho \, (1-\alpha)$ 
in~(\ref{eqBolmua}) is fundamental in order to get uniformity of these bounds in the limit $\alpha \to 1$. 
Let us also mention that our choice of scaling for 
the equation~(\ref{eqBolmua}) is mass invariant, that is $F$ with density $\rho(F)$ satisfies the equation 
if and only if $F/\rho(F)$ satisfies the equation with $\rho=1$. Therefore all the estimates on the profiles 
are homogeneous in terms of the density $\rho$.  

  \begin{prop}\label{estimatesonGe} 
  Let us fix $\alpha_0 \in (0,1)$. There exists $a_1,a_2,a_3,a_4 \in (0,\infty)$ and, for any $k \in \N$,
  there exists $C_k \in (0,\infty)$ such that 
\beqn\label{estimGeunif}\qquad
\forall \, \alpha \in [\alpha_0,1), \,\, \forall \, F_\alpha \in \FF_\alpha, \quad
\| F_\alpha \|_{L^1(  e^{a_1 \, |v| })} \le a_2,  \quad  \| F_\alpha \|_ {H^k(\RR^N)} \le C_k.
\eeqn
\end{prop}


\medskip\noindent{\sl Proof of Proposition \ref{estimatesonGe}.} 
We split the proof into several steps. We fix $\alpha \in [\alpha_0,1)$ 
and $F_\alpha \in \FF_\alpha$ for which we will establish the announced bounds. 
Let emphasize that thanks to the {\it a priori} bounds satisfied by $F_\alpha$ all the computations 
we will perform are rigorously justified.  From now we omit the subscript ``$\alpha$'' when no confusion is possible. 

\medskip\noindent
{\it Step~1. Upper bound on the energy using the energy dissipation term. } 
We prove that 
\beqn\label{UpperEnergy}
\forall \, \alpha \in (0,1] \qquad \EE \le \rho \, \left( {2 \, N \over b_1} \right)^{2/3}. 
\eeqn
From equation~(\ref{EnergyDEGe}) on the energy of the profile $G$ there holds
\beqn\label{dissipenergyGe}
 (1+\alpha) \, b_1  \, \int_{\RR^N} \! \int_{\RR^N} F \, F_* \, |u|^3 \, dv \, dv_* = 
 2 \, N \, \rho^2. 
\eeqn
From Jensen's inequality
  $$
  \int_{\R^N} |u|^3 \, F_* \, dv_* \ge \rho \, |v|^3,
  $$
and H\"older's inequality
  $$
  \int_{\R^N} |v|^3 \, F \, dv \ge \rho^{-1/2} \, \left( \int_{\R^N} |v|^2 \, F \, dv \right)^{3/2}, 
  $$
we get 
  $$
  (1+\alpha) \, b_1 \, \rho^{1/2} \, \EE^{3/2} \le 2 \, N \, \rho^2
  $$ 
from which the bound (\ref{UpperEnergy}) follows. 

\medskip\noindent
{\it Step 2. Lower bound on the energy using the entropy.} We prove 
\beqn\label{LowerEnergy}
\forall \, \alpha \in (0,1] \qquad \EE \ge \rho \, \left( {\alpha^2 \, N^2 \over \sqrt{2} \, b_2 } \right)^{2/3}  
\quad\hbox{with}\quad
b_2 := \| b \|_{L^1}.
\eeqn
  
\begin{rem}
The choice of scaling we have made for the rescaled equation~(\ref{eqBolmua})  
becomes clear from this computation: it is chosen 
such that the energy of the stationary solution does not blow up nor vanishes for $\alpha \to 1$. 
The restriction $\alpha \in [\alpha_0,1)$, $\alpha_0 > 0$, is then made in order to get a uniform 
estimate from below on the energy. 
\end{rem}

\smallskip
By integrating the equation satisfied by $F$ against $\log F$ we find   
  $$
  \int_{\R^N} Q(F,F) \log F \, dv 
+ \rho \, (1-\alpha) \,  \int_{\R^N} \log F  \,\,  \Delta_v F \, dv = 0.
  $$
Then we write the first term as in~\cite[Section~1.4]{GPV**} to find 
  \bean
  &&{1 \over 2} \int\!\! \int\!\! \int_{\R^{2N}\times S^{N-1}} F \, F_* \left( \log {F'F'_* \over F F_*} 
  - {F'F'_* \over F F_*} + 1 \right) \, B \, dv \, dv_* \, d\sigma \\
  &&\qquad + {1 \over 2} \int\!\! \int\!\! \int_{\R^{2N}\times S^{N-1}}   \left( F'F'_* -F F_* \right) \, B \, dv \, dv_* \, d\sigma  
  - \rho \,(1-\alpha) \, \int_{\R^N} {|\nabla_v F |^2 \over F} \, dv = 0. 
  \eean
Recalling that $x - \log \, x - 1 \ge 0$ for any $x \ge 0$ and making the change of variables 
$(v',v'_*) \to (v,v_*)$ in the second term we obtain
\beqn\label{DissipGe} 
\rho \,(1-\alpha) \, \int_{\R^N} {|\nabla_v F |^2 \over F} \, dv  
\le  {1 \over 2} \left(  {1 \over \alpha^2} - 1 \right) \, b_2 \, \int \!\! \int_{\R^{2N}}  F \, F_* \, |u| \, dv \, dv_*.
\eeqn
On the one hand, from Cauchy-Schwarz's inequality
\bear\label{DissipGe2} \nonumber
&& \int \!\! \int_{\R^{2N}} F \, F_* \, |u| \, dv \, dv_* \le \\
&& \le \left( \int \!\!  \int_{\R^{2N}} F \, F_*  \, dv \, dv_* \right)^{1/2} \!\! 
 \left( \int \!\!  \int_{\R^{2N}} F \, F_* |u|^2 \, dv \, dv_* \right)^{1/2} = \sqrt{2}  \, \rho^{3/2}\, \EE^{1/2}.
\eear
On the other hand, we compute 
\bear\label{DissipGe3} \nonumber
&& 0 \le \int_{\R^N} \left| 2 \, \nabla \sqrt{F} + {\rho \, N \, v \over \EE} \, \sqrt{F}\right|^2 \, dv  \\ \nonumber
&& \quad = \int_{\R^N} \left(4 \, {|\nabla \sqrt{F}|^2 }+ 2 \,  {\rho \, N \, v \over \EE} \cdot \nabla F 
  + \left( {\rho \, N  \over \EE}\right)^2 \, |v|^2 \, F\right)\, dv \\
&& \quad \le  2 \, \int_{\R^N} {|\nabla_v F |^2 \over F} \, dv  - 2 \,  {(\rho \, N)^2  \over \EE} 
                  + \left( {\rho \, N  \over \EE}\right)^2 \, \EE= 
 2 \, \int_{\R^N} {|\nabla_v F |^2 \over F} \, dv -  {(\rho \, N)^2  \over \EE}.
\eear
Gathering (\ref{DissipGe}), (\ref{DissipGe2}) and (\ref{DissipGe3}) we get
$$
\rho \, {(\rho \, N)^2  \over 2 \, \EE} \le {1 \over 2}\, {1+\alpha \over \alpha^2} \, b_2  \sqrt{2}  \, \rho^{3/2}\, \EE^{1/2},
$$
from which we deduce   (\ref{LowerEnergy}). 

\medskip\noindent
{\it Step 3. Upper bound on exponential moments.} 
There exists $A,C > 0$ such that 
  $$ 
  \forall \, \alpha \in [0,1), \qquad \int_{\R^N}  F (v) \, e^{A |v|} \, dv \le C \, \rho.
  $$

We refer to \cite[Theorem 1]{BGP**} where this bound is obtained as an immediate consequence 
of the following sharp moment estimates: there exists an explicit $X > 0$ such that 
  \beqn\label{Uppermoment}
  \forall \, \alpha \in [0,1), \qquad {\bf m}_k = \int_{\R^N} G \, |v|^k \, dv \le \Gamma(k+1/2) \, X^{k/2} \, \rho.
  \eeqn
It is worth noticing that in \cite{BGP**} the Povzner inequality used 
in order to get (\ref{Uppermoment}) is uniform in terms of the normal restitution  
coefficient $\alpha \in [0,1]$ appearing in the collisional operator and in terms of the 
coefficient $\mu$ in front of the thermal bath term. Moreover the factor $\rho$ comes from our 
choice of the scaling variables (in which $\rho$ is involved). 

\medskip\noindent
{\it Step~4.  Uniform upper bound on the $L^2$ norm, and $H^k$ norms, $k >0$.} 
On the basis of the uniform bounds from below and above on the energy, 
the proof can be done exactly as in \cite[Proposition~2.1]{MMGranular3}, 
since the diffusion term plays no role in the energy estimates (it only helps). \qed

\medskip
After the estimates from above on the profile, let us now state a pointwise bound 
from below (it is a straightforward consequence of Proposition~\ref{LowerBdga} in 
the appendix).

\begin{prop}\label{estimatesonGe2} 
Let us fix $\alpha_0 \in (0,1)$. There exists $a_3,a_4 \in (0,\infty)$ such that 
\beqn\label{estimGeunif2}\qquad
\forall \, \alpha \in [\alpha_0,1), \,\, \forall \, F_\alpha \in \FF_\alpha, \quad
F_\alpha \ge a_3 \, e^{-a_4 \, |v|^{8}}.
\eeqn
\end{prop}    

\subsection{Estimates on the difference of two stationary solutions}

In this subsection we take advantage of the mixing effects of the collision 
operator in order to show that the $L^1$ norm of their difference of two stationary solutions 
(corresponding to the same inelasticity coefficient) indeed controls  
the $H^k \cap L^1(m^{-1})$ norm of their difference for any $k \in \N$ and for 
some exponential weight function $m$, uniformly in terms of $\alpha \in [\alpha_0,1)$. 

 \begin{prop}\label{errorH1} 
 For any $k >0$, there is $m=\exp(-a \, |v|)$, $a \in (0,\infty)$ and 
 $C_{k}>0$ such that for any $\alpha \in [\alpha_0,1)$ 
 and any $F_\alpha,H_\alpha \in \FF_\alpha$ there holds
   \beqn\label{eq:errorH1}
   \|H_\alpha - F_\alpha  \|_{H^k \cap L^1(m^{-1})} \le C_{k} \, \|H_\alpha - F_\alpha \|_{L^1}. 
   \eeqn
 \end{prop}

\smallskip\noindent{\sl Proof of Proposition \ref{errorH1}. }
We proceed in three steps. It is worth mentioning that all the constants in the proof are uniform in terms of the 
normal restitution coefficient $\alpha \in [\alpha_0,1)$, as they only depend on the 
uniform bounds of Proposition~\ref{estimatesonGe} and some uniform bounds on the 
collision kernel. 

\medskip\noindent
The proof is based on the following three steps, which can all be proved following 
exactly the arguments in~\cite[Proposition~2.7]{MMGranular3}. 
 
{\sl Step 1. Control of the $L^1$ moments.} 
We prove first that there exists $A, C \in (0,\infty)$ such that 
  $$
  \forall \, \alpha \in [\alpha_0,1), \qquad  \int_{\R^N} |H_\alpha-F_\alpha| \, e^{A \, |v|} \, dv 
      \le C  \int_{\R^N} |H_\alpha-F_\alpha| \, dv.
  $$
The proof can be done exactly as in the Step~1 of~\cite[Proposition~2.7]{MMGranular3}: the 
only change is the replacement of the vanishing anti-drift term by the term viscosity term, which 
both in any case are negligible in the moments estimates (note also that this estimate on the 
diffusive Boltzmann equation is proved within the paper~\cite{BGP**}, with constant 
independent on the viscosity coefficient).

\medskip\noindent
{\sl Step 2. Control of the $L^2$ norms.} 
It is done exactly as in Step~2 of~\cite[Proposition~2.7]{MMGranular3}, using the 
regularity theory of the gain term: the only difference is that instead anti-drift term 
which vanishes in the estimates, one has now a diffusion term which yields good damping 
negative in the $L^2$ energy estimate (which can therefore be dropped). 

\medskip\noindent
{\sl Step 3. Control of the $H^k$ norms.} 
The proof can be done exactly as in Step~3 of~\cite[Proposition~2.7]{MMGranular3}. 
The anti-drift is replaced by the diffusion term which is also a good negative term 
in the $H^k$ energy estimates.
\qed


\section{Quantification of the elastic limit $\alpha \to1$}
\setcounter{equation}{0}
\setcounter{theo}{0}


We have the following estimate on the distance between $F_\alpha$ and $\bar F_1$ for any stationary solution $F_\alpha$. 
\smallskip
  \begin{prop}\label{ConvGexplicit} 
  For any $\eps > 0$ there exists $C_\eps$ (independent of the mass $\rho$) such that 
    \beqn\label{estimateonGe2}
    \forall \, \alpha \in [\alpha_0,1) \qquad  \sup_{F_\alpha \in \FF_\alpha} 
    \, \|F_\alpha - \bar F_1 \|_{L^1_2} \le C_\eps \, \rho \, (1-\alpha)^{{1 \over 2+\eps}}
    \eeqn
where we recall that  $\bar F_1$ is the Maxwellian function defined by (\ref{Maxlim})--(\ref{tempSS}).
  \end{prop}

\medskip
\noindent{\sl Proof of Proposition \ref{ConvGexplicit}.} 
On the one hand, for any inelasticity coefficient $\alpha \in [\alpha_0,1)$ and 
profile $F_\alpha$, there holds from (\ref{DissipGe}) together 
with Corollary~\ref{DHeDH1} and the uniform estimates of Proposition~\ref{estimatesonGe}
\beqn\label{DH1<DHa}
D_{H,1}(F_\alpha)  \le D_{H,\alpha}(F_\alpha) + \rho^2 \, \OO(1-\alpha) \le \rho^2 \, \OO(1-\alpha).
\eeqn
On the other hand, introducing the Maxwellian function $M_\theta$ with the same mass, momentum and temperature 
as $F_\alpha$, that is $M_\theta$ given by (\ref{defMrut}) with $u=0$ and 
$\theta = \EE(F_\alpha)/\rho$, and gathering (\ref{DH1<DHa}), (\ref{EEPineg}), (\ref{CKineg}) 
with the uniform estimates of Proposition~\ref{estimatesonGe} and interpolation inequality, 
we obtain that for any $q,\eps > 0$ there exists $C_{q,\eps}$ such that 
\beqn\label{estimGeMtheta}
\forall \alpha \in [\alpha_0,1) \qquad \| F_\alpha - M_\theta \|_{L^1_q}^{2+\eps} \le C_{q,\eps} \, 
\rho^{2+\eps} \, (1-\alpha).
\eeqn
Next, from (\ref{dissipenergyGe}), we  have 
$$
 b_1  \, \int_{\RR^N} \! \int_{\RR^N} F_\alpha \, F_{\alpha*} \, |u|^3 \, dv \, dv_* - 2 \, \rho \, N \, 
\int_{\R^N} F_\alpha \,  dv = 
 (1-\alpha)  \,  {b_1 \over 2} \, \int_{\RR^N} \! \int_{\RR^N} F_\alpha \, F_{\alpha*} \, |u|^3 \, dv \, dv_*
$$
and then 
\beqn\label{ineqtheta1}
\left|  \Psi(\theta) \right| \le C_1 \, \| F_\alpha - M_\theta \|_{L^1_3} + C_2 \, \rho^2 \, (1-\alpha),
\eeqn
where we have used that $F_\alpha$ and $M_\theta$ are bounded thanks to Proposition~\ref{estimatesonGe} and we have defined 
\beqn\label{defPsi}
\Psi(\theta) = 2 \, \rho \, N \, \int_{\R^N} M_\theta \,  dv  
- b_1 \, \int_{\RR^N} \! \int_{\RR^N} M_\theta \, M_{\theta*} \, |u|^3 \, dv \, dv_*.
\eeqn
By elementary changes of variables, this formula simplifies into 
$$
\Psi(\theta) = k_1 - k_2 \, \theta^{3/2}
$$
with $k_1 = 2 \, \rho^2 \, N$ and, using (\ref{MMu3}),
$$
k_2 =   \rho^2 \, b_1 \, \int_{\RR^N \times \RR^N} M_{1,0,1} \, (M_{1,0,1}) _* \, |u|^3 \, dv \, dv_* =
2^{3/2} \, \rho^2 \, b_1 \, {\bf m}_{3/2}(M_{1,0,1}).
$$ 
We next observe that $\Psi \in C^\infty(0,\infty)$ and $\Psi$ is strictly concave. 
It is also obvious that the equation $\Psi(\theta) = 0$ for $\theta > 0$ has a unique solution 
which is $\bar\theta_1$ defined in~(\ref{tempSS}), and that we have 
$$
\Psi(\theta) \le \Psi'(\bar\theta_1) \, (\theta-\bar\theta_1) = - C \, (\theta-\bar\theta_1)/2
$$ 
for some explicit constant, as well as
\beqn\label{propPsi}
\Psi(\theta) = k_2 \, [\bar\theta_1^{3/2} -  \theta^{3/2}].
\eeqn
Plugging this expression for $\Psi$ into~(\ref{ineqtheta1}) and using the lower bound (\ref{LowerEnergy}) 
on the temperature $\theta$ and the estimate (\ref{estimGeMtheta}) we obtain that for any $\eps > 0$ 
there is $C_\eps \in (0,\infty)$ such that 
\beqn\label{thetaa-b}
\forall \,  \alpha \in (\alpha_0,1) \qquad 
 \left| \theta^{3/2} -  \bar\theta_1^{3/2} \right| ^{2+\eps} \le C_\eps \, (1-\alpha). 
\eeqn
Namely, we have thus proved that the temperature of $\bar F_\alpha$ converge (with rate) to the 
expected temperature $\bar\theta_1$. In order to come back to the norm of $F_\alpha - \bar F_1$, 
we first write, using Cauchy-Schwarz's inequality,
\bear \nonumber
\| F_\alpha - \bar F_1 \|_{L^1_{-N}} 
&\le& \| F_\alpha - M_\theta \|_{L^1_{-N}} +  \| M_\theta - \bar F_1 \|_{L^1_{-N}} \\ \label{Ga-G1L1}
&\le&  \| F_\alpha - M_\theta \|_{L^1} + C_N \, \| M_\theta - \bar F_1 \|_{L^2}, 
\eear
and we remark that 
\beqn\label{Mtheta-G1L2}
\| M_\theta - \bar F_1 \|_{L^2}^2 
\le C \, \rho^2 \, \left|\theta^{3/2} - \bar\theta_1^{3/2}\right|. 
\eeqn
Gathering (\ref{Ga-G1L1}) with  (\ref{Mtheta-G1L2}), (\ref{thetaa-b}) and (\ref{estimGeMtheta}) we deduce 
that for any $\eps > 0$ there is $C_\eps \in (0,\infty)$ such that 
\bean
\forall \,  \alpha \in (\alpha_0,1) \qquad  
\| F_\alpha - \bar F_1 \|_{L^1_{-N}}^{2+\eps} \le C_\eps \, \rho^{2+\eps} \, (1-\alpha), 
\eean
and (\ref{estimateonGe2}) follows by interpolation again. 
\qed
 

\section{Uniqueness and continuity of the path of stationary solutions}
\setcounter{equation}{0}
\setcounter{theo}{0}

\begin{theo}\label{uniqueness1} 
There exists a constructive $\alpha_1 \in (0,1)$ 
such that the solution $F_\alpha$ of (\ref{eqBolmua}) is unique for any $\alpha \in [\alpha_1,1]$. 
We denote by $\bar F_\alpha$ this unique stationary solution. 
\end{theo}

That is an immediate consequence of the following result. 

\begin{prop}\label{uniqueness2} 
There is a constructive constant $\eta \in (0,1)$ such that 
$$ 
\quad \left.
  \begin{array}{l}
  G, \, H \in \FF_\alpha, \,\, \alpha \in (1-\eta,1) \vspace{0.3cm} \\
  \|G - \bar F_1 \|_{L^1_2} \le \eta, \,\,\,
    \|H - \bar F_1 \|_{L^1_2} \le \eta 
    \end{array}
  \right\} \quad\hbox{implies}\quad G = H. 
$$
\end{prop}

\smallskip\noindent{\sl Proof of Theorem~\ref{uniqueness1}. }
Let us assume that Proposition~\ref{uniqueness2} holds. 
Then Proposition~\ref{ConvGexplicit} implies that there is 
some explicit $\e \in (0,1)$ such that for $\alpha \in (1-\e,1]$ 
one has 
$$
\sup_{F_\alpha \in \FF_\alpha} \| F_\alpha - \bar F_1 \|_{L^1 _2} \le \eta
$$
where $\eta$ is defined in the statement of Proposition~\ref{uniqueness2}. Up to reducing 
$\eta$, it is always possible to take $\eta \le \e$, and the proof 
is completed by applying Proposition~\ref{uniqueness2}. 
\qed

\smallskip\noindent{\sl Proof of Proposition~\ref{uniqueness2}. }
Let us consider any exponential weight function  $m$ with $s \in (0,1)$, $a \in (0,+\infty)$, 
or with $s=1$ and $a \in (0,\infty)$ small enough. 
Let us also define the subvector space of $\Ll^1(m^{-1})$ of functions with zero energy 
\beqn\label{defOO}
\OO= \Cc_{0,0,0}  \cap \Ll^1(m^{-1}),
\eeqn
the function $\psi = C \, (|v|^2 - N) \, M_{1,0,1}$ 
such that its mass is zero and its energy is $\EE(\psi) = 1$, and $\Pi$ the following projection 
$$
\Pi : \Ll^1(m^{-1}) \to \OO, \qquad \Pi(g) = g - \EE(g) \, \psi.
$$
Finally, let us introduce $\Phi$ the following non-linear functional operator 
$$
\Phi : [0,1) \times ( W^{2,1}_1(m^{-1}) \cap \Cc_{\rho,0} ) \,\, \to \,\,  \R \times  \OO,
$$
and
$$
\Phi(1, \cdot ) : ( L^1_1(m^{-1}) \cap \Cc_{\rho,0} ) \,\, \to \,\,  \R \times  \OO,
$$ 
by setting 
\bean
\Phi(\alpha,g) = \left(  (1+\alpha) D_\EE(g)  -  2 \, N \,  \rho^2 , \,  \Pi \Bigl[Q_\alpha(g,g) 
        + \tau_\alpha \, \Delta_v g \Bigr]  \right). &&
\eean
It is straightforward that $\Phi(\alpha,F_\alpha) =0$ for any $\alpha \in [\alpha_0,1]$ and 
$F_\alpha \in \FF_\alpha$, and that the equation 
$$
\Phi(1,g) = (0,0)
$$
has a unique solution, given by $g = \bar F_1 = M_{\rho,0,\bar\theta_1}$ defined in (\ref{Maxlim}), (\ref{tempSS}). 

The function $\Phi$ is linear and quadratic in its second argument by inspection, and easy computations 
yield the following formal differential according to the second argument at the point $(1,\bar F_1)$: 
\beqn\label{D2Phi}
  D_2 \Phi(1,\bar F_1) \, h =: A \, h := \Bigg(  4 \,  \tilde D_\EE( \bar F_1, h),  \ 2 \, \tilde Q_1(\bar F_1,h) \Bigg) 
\eeqn
where $\tilde Q_\alpha$ is defined in (\ref{Qinelsym}) and 
$$
\tilde D_\EE (g,h) := {b_1} \, \int\!\!\int_{\R^N \times \R^N} g \, h_*\, |u|^3 \, dv \, dv_*.
$$
Notice that we can remove the projection on the last argument in (\ref{D2Phi}) since the elastic collision 
operator always has zero energy. 

\medskip\noindent
On the basis of Lemma~\ref{lemDPhiM} in the appendix for $A$, we shall now 
prove Proposition~\ref{uniqueness2}. We write 
\bear
F_\alpha - H_\alpha &=& A^{-1} \, \Big[ A \, F_\alpha - \Phi(\alpha,F_\alpha) 
      + \Phi(\alpha,H_\alpha) - A \, H_\alpha  \Big]  \nonumber \\
\label{DA-1I}
&=& A^{-1} \, (I_1, I_2) 
\eear
with (recall that the bilinear operators $\tilde D_\EE$ and $\tilde Q_\alpha$ are symmetric)
\[
\left\{
\begin{array}{l}
I_1 :=  4 \, \tilde D_\EE(\bar F_1, F_\alpha - H_\alpha) 
- (1+\alpha) \, D (F_\alpha) + (1+\alpha) \, D (H_\alpha) \vspace{0.3cm} \\
I_2 := \Pi \, I_{2,1}+  \Pi \, I_{2,2}
\end{array}
\right.
\]
and
\[
\left\{
\begin{array}{l}
I_{2,1} :=  2 \, \tilde Q_1(\bar F_1,F_\alpha-H_\alpha) 
                - Q_\alpha(F_\alpha,F_\alpha) + Q_\alpha(H_\alpha,H_\alpha)  \vspace{0.3cm} \\
I_{2,2} :=  \rho \, (1-\alpha)  \, \Delta (H_\alpha - F_\alpha). 
\end{array}
\right.
\]
On the one hand, 
$$
I_1 = 2 \, D\big(2 \bar F_1 - (F_\alpha + H_\alpha), F_\alpha - H_\alpha\big) 
+ (1-\alpha)  \, D (F_\alpha + H_\alpha, F_\alpha - H_\alpha)
$$
so that 
\bear \nonumber
|I_1| &\le& C_3 \, 
\Big( \|\bar F_1 - F_\alpha\|_{L^1 _3} + \|\bar F_1 - H_\alpha \| _{L^1 _3} \\ \nonumber
&& \hspace{4cm} + (1-\alpha) \, 
 \|F_\alpha \|_{L^1 _3} + (1-\alpha) \, \|H_\alpha \|_{L^1 _3} \Big) \, \|F_\alpha - H_\alpha \|_{L^1 _3}
\nonumber \\
\label{I2} &\le& \eta_1(\alpha) \, \|F_\alpha - H_\alpha \|_{L^1_1(m^{-1})}
\eear
with $\eta_1(\alpha) \to 0$ when $\alpha \to1$ (with explicit rate, for instance $\eta_1(\alpha) = C_1 \, (1-\alpha)^{1/3}$) 
because of Propositions~\ref{estimatesonGe} and~\ref{ConvGexplicit}.

On the other hand, 
\bean
I_{2,1} &=& 2 \, \tilde Q_1(\bar F_1, F_\alpha - H_\alpha) -  2 \, \tilde Q_\alpha (\bar F_1, F_\alpha - H_\alpha) 
+ Q_\alpha(\bar F_1 - F_\alpha, F_\alpha - H_\alpha) + Q_\alpha (F_\alpha - H_\alpha,\bar F_1 - H_\alpha). 
\eean  
From Proposition~\ref{estimatesonGe}  and estimate (\ref{convLetoL1}) in  Proposition~\ref{ContQe1} there holds 
  \[ \| \tilde Q_1(\bar F_1, F_\alpha - H_\alpha) -  \tilde Q_\alpha(\bar F_1, F_\alpha - H_\alpha) \|_{L^1(m^{-1})} \le 
        C\, (1 - \alpha) \, \| F_\alpha - H_\alpha \|_{W^{1,3}_3(m^{-1})}. \]
From   estimate (\ref{QaWk1}) in  Proposition~\ref{ContQe1}  we have 
\bean 
&& \|  Q_\alpha(\bar F_1 - F_\alpha, F_\alpha - H_\alpha) 
+ Q_\alpha (F_\alpha - H_\alpha,\bar F_1 - H_\alpha) \|_{L^1(m^{-1})} \\ 
&& \hspace{3cm} 
\le C_4 \, \left( \| F_\alpha - \bar F_1 \|_{L^1 _1(m^{-1})} + \| H_\alpha -\bar F_1 \|_{L^1 _1(m^{-1})} \right) 
           \, \| F_\alpha - H_\alpha \|_{L^1 _1(m^{-1})}.
\eean
Together with Propositions \ref{ConvGexplicit} we thus obtain
\beqn \label{I31}
\|I_{2,1} \|_{L^1(m^{-1})} \le \eta(\alpha) \, \|F_\alpha - H_\alpha \|_{L^1_1(m^{-1})}
\eeqn 
for some $\eta(\alpha) \to 0$ as $\alpha \to 1$. Here we can take for instance 
(when $s=1/2$ in the formula of $m$)
$\eta(\alpha) = C \, (1-\alpha)^{1/12}$ 
for some $C \in (0,\infty)$ by picking a suitable $\varepsilon$ and interpolating. 

Finally from Proposition~\ref{errorH1} there holds 
\beqn \label{I32}
\|I_{2,2} \|_{L^1(m^{-1})} \le C_5 \, (1-\alpha) \, \|F_\alpha - H_\alpha \|_{L^1_1(m^{-1})}.
\eeqn 
Gathering~(\ref{I2}), (\ref{I31}) and (\ref{I32}) we obtain from (\ref{DA-1I}), Lemma~\ref{lemDPhiM} and Proposition~\ref{errorH1} again
$$
\| F_\alpha - H_\alpha \|_{L^1_1(m^{-1})} \le \eta(\alpha) \, \|A^{-1}\| \, \|F_\alpha - H_\alpha \|_{L^1_1(m^{-1})}
$$
for some function $\eta$ such that  $\eta(\alpha) \to 0$ as $\alpha \to 1$ (with explicit rate). 
Hence choosing $\alpha_1$ close enough to $1$ 
we have $\eta(\alpha) \, \|A^{-1}\|  \le 1/2$ for any $\alpha \in [\alpha_1,1)$. 
This implies $F_\alpha = H_\alpha$ and concludes the proof. \qed 

\medskip
Let us state now without proof, because is is completely similar to the corresponding result in ~\cite{MMGranular3} and is again a variation around the implicit function theorem, some  regularity result on the path of self-similar profiles.

\begin{lem}\label{Contderivee=1Ge}  
The map $[\alpha_1,1] \to L^1(m^{-1})$, $\alpha \mapsto \bar F_\alpha$ is continuous 
on $[\alpha_1,1]$ and differentiable at $\alpha= 1$. More precisely, 
there exists $\bar F'_1 \in L^1(m^{-1})$ and for any $\eta \in (1,2)$ 
there exists a constructive $C_\eta \in (0,\infty)$ such that 
\beqn\label{devGalpha1}
\| \bar F_\alpha - \bar F_1 - (1-\alpha) \, \bar F'_1 \|_{L^1(m^{-1})} \le C_\eta \, (1-\alpha)^\eta
\qquad \forall \, \alpha \in (\alpha_0,1).
\eeqn
\end{lem}


\section{Spectral study of the linearized problem} \label{sec:spectral}
\setcounter{equation}{0}
\setcounter{theo}{0}

In this section we shall enounce some results on the geometry of the spectrum of the 
linearized diffusive inelastic collision operator 
for a small inelasticity, as well as  estimates on its resolvent and on the associated linear semigroup. 
This is straightforwardly adapted from \cite{MMGranular3}.

We thus consider the operator 
  \[ \LL_\alpha : f \mapsto Q_\alpha(f,f) + \tau_\alpha \, \Delta f \]
and some fluctuations $h$ around the stationary solution $\bar F_\alpha$:  
that means $f = \bar F_\alpha  +  h$ with $f \in L^1(m^{-1})$ where $m$ is a fixed smooth exponential weight function, 
as defined in (\ref{defdem}).  
The corresponding linearized unbounded operator $\LL_\alpha$ is acting on $L^1(m^{-1})$  
with domain $\hbox{dom}(\LL_\alpha) = W^{2,1}_1(m^{-1})$ if $\alpha \not= 1$ and 
$\hbox{dom}(\LL_1) = L^1_1(m^{-1})$ (it is straightforward to check that it is closed in this space).
Since the equation in self-similar variables preserves mass and the zero momentum, the correct spectral 
study of $\LL_\alpha$ requires to restrict this operator to zero mean and centered distributions 
(which are preserved as well), that means to work in $\Ll^1(m^{-1})$. 
When restricted to this space, the operator $\LL_\alpha$ is denoted 
by $\hat \LL_\alpha$. We denote by $R(\hat \LL_\alpha)$ the resolvent 
set of $\hat \LL_\alpha$, and by $\Rr_\alpha(\xi) = (\hat \LL_\alpha - \xi)^{-1}$ its 
resolvent operator for any $\xi \in R(\hat \LL_\alpha)$. 

\subsection{The result on the spectrum and resolvent}
\label{subsec:spec}

The result proved in this section is a translation of  
\cite[Theorem~5.2]{MMGranular3} for the diffusive inelastic Boltzmann 
equation. Let us define for any $x \in \R$ 
the half-plane $\Delta_x$ by
$$
\Delta_x = \{ \xi \in \CC, \,\, \Re e \, \xi \ge x \}.
$$

\begin{theo}\label{thLalpha} 
Let us fix $\bar \mu \in (\mu_2,0)$, with $\mu_2 < 0$ defined in Theorem~\ref{theoL1}, $k,q \in \N$ and $m$ a smooth weight  exponential function with $s \in (0,1)$. 
Then there exists $\alpha_2 \in (\alpha_1,1)$ such that for any 
$\alpha \in [\alpha_2,1]$ the following holds:

\begin{itemize}
\item[(i)] The unbounded operator $\hat\LL_\alpha$ is well defined and closed in $\Ww^{k,1} _q(m^{-1})$ ($\forall \, k \ge 0$, $\forall \, q \ge 0$). Its spectrum $\Sigma(\hat  \LL_\alpha)$ satisfies
$$
\Sigma(\hat \LL_\alpha) \subset \Delta_{\bar \mu} ^c \cup \{ \mu_\alpha \},
$$
where $\mu_\alpha$ is a real eigenvalue which does not depend on the choice of the space $\Ww^{k,1} _q(m^{-1})$  and satisfies (with explicit bounds)
\beqn\label{expanmua}
\mu _\alpha = - 3 \, \rho \, (1-\alpha) + \OO(1-\alpha)^2 \quad\hbox{when}\quad \alpha \to 1.
\eeqn
Moreover,  $\mu_\alpha$  is a $1$-dimensional eigenvalue, and more precisely, the 
eigenspace associated to $\mu_\alpha$ is $\R \, \phi_\alpha$ with 
$\| \phi_\alpha \|_{L^1_2} = 1$ and 
$\| \phi_\alpha -  \phi_1 \|_{\Ww^{k,1} _q(m^{-1})} \le C \, (1-\alpha)$, 
where $\phi_1 := c_0 \, (|v|^2 - N \, \bar\theta_1) \, \bar F_1$ is the ``energy eigenfunction'' 
associated to the linearized  elastic Boltzmann operator.

\item[(ii)] The resolvent  $\Rr_\alpha(\xi)$  in $\Ww^{k,1}_q(m^{-1})$ is holomorphic on a 
neighborhood of $\Delta_{\bar \mu} \backslash \{ \mu_\alpha \}$ and 
there are explicit constants $C_1,C_2$ such that 
$$
\sup_{z \in \CC, \,\, \Re e \, z = \bar \mu} 
\| \Rr_\alpha(z)|\|_{\Ww^{k,1} _q(m^{-1}) \to \Ww^{k,1} _q(m^{-1})} \le C_1
$$
and 
$$
\| \Rr_\alpha(\bar \mu + i s) \|_{\Ww^{k+2,1} _{q+1}(m^{-1}) \to \Ww^{k,1} _q(m^{-1})} \le 
\frac{C_2}{1+|s|}.
$$

\item[(iii)] The linear semigroup $S_\alpha(t)$ associated to $\hat \LL_\alpha$ in 
$\Ww^{k,1} _q (m^{-1})$ writes
$$
S_\alpha(t) = e^{\mu_\alpha \, t}  \, \Pi_\alpha + R_\alpha(t),
$$
where $\Pi_\alpha$ is the projection on the ($1$-dimensional) eigenspace associated to $\mu_\alpha$ 
and where $R_\alpha(t)$ is a semigroup which satisfies 
\beqn\label{estimResta}
\|R_\alpha(t) \|_{\Ww^{k+2,1} _{q+2}(m^{-1}) \to \Ww^{k,1} _q(m^{-1})}
 \le C_{k} \, e^{\bar \mu \, t } 
\eeqn
with explicit bounds. 
\end{itemize}
\end{theo}

\subsection{Sketch of proof}
\label{subsec:proofspec}

The proof is straightforwardly adapted from \cite[Section~5]{MMGranular3}. 
We shall only mention the main steps of the proof and we emphasize the few points which differs here (due to the replacement  of the anti-drift term by a diffusive term). 

The proof is based on results on the elastic collision operator together with the 
decomposition for $\xi  \in \mathbb C$
\beqn\label{decomp1}
\LL_\alpha - \xi = A_\delta - B_{\alpha,\delta} 
\eeqn
where 
\beqn\label{decomp2}
A_\delta = \LL^+_{1,\delta} - \LL^*, \quad
B_{\alpha,\delta} = \nu + \xi + \big( \LL^+_{1,\delta} - \LL_1 ^+ \big) + P_\alpha
\eeqn
with $P_\alpha = \LL_1 - \LL_\alpha = \LL_1^+ - \LL_\alpha^+ + \tau_\alpha \, \Delta$, 
and where $\LL_1 = \LL^+ - \LL^* -\nu$ is the usual decomposition of the elastic 
collision operator and $\LL^+ _{1,\delta}$ is the 
regularized truncation of the ``gain'' part introduced in \cite{GM:04}.

The differences with the similar decomposition in \cite{MMGranular3} only lies in the $B_{\alpha,\delta}$ and $P_\alpha$ opertors.  But the key technical estimates are still true for this operator. 
Indeed it can be checked straightforwardly that the following holds (let us point out that 
one loses now two indices of regularity instead  of one, which explains the minor change 
in the statement of the estimate on the resolvent). 

\begin{lem}[See Lemmas~5.8 \& 5.9 in \cite{MMGranular3}]\label{techLaBdelta} 
Let us fix $k,q \ge 0$ and an exponential weight function $m$. 

\begin{itemize}
\item[(i)] There exists some constant $C$ such that for any $\alpha \in (\alpha_0,1]$
$$
\left\|  \LL_\alpha \right\|_{W^{k+2,1}_{q+1} (m^{-1}) \to W^{k,1}_q(m^{-1})} \le C, \quad
\left\| \LL_1 - \LL _\alpha \right\|_{W^{3,1}_{3} (m^{-1}) \to L^1(m^{-1})}  \le C \, (1- \alpha).
$$

\item[(ii)] There exists some constants $\delta^* > 0$ and  $\alpha_2 \in (\alpha_1,1)$ such that for any $\xi \in \Delta_{\mu_2}$,  $\delta \in [0,\delta^*]$ and $\alpha \in [\alpha_2,1]$, the operator 
$$
B_{\alpha,\delta} : W^{k+2,1} _{q+1}(m^{-1}) \to W^{k,1}_q (m^{-1})
$$ 
is invertible and the  inverse operator $B_{\alpha,\delta}(\xi) ^{-1}$ satisfies 
$$
\left\| B_{\alpha,\delta}(\xi) ^{-1} \right\|_{W^{k,1}_q (m^{-1}) \to W^{k,1}_q (m^{-1})} 
\le \frac{C_1}{\mbox{{\em dist}}(\Re e \, \xi,\nu(\R^N))} 
$$
and 
$$
\left\| B_{\alpha,\delta}(\xi) ^{-1} \right\|_{W^{k,1}_{q} (m^{-1}) \to W^{k+2,1}_{q}(m^{-1})}
\le \frac{C_2}{\mbox{{\em dist}}(\xi,\nu(\R^N))}
$$
for some explicit constants $C_1,C_2>0$ depending on $k,q,\delta^*,\alpha_2$.

\item[(iii)] As a consequence, the resolvent operator $\Rr_\alpha(\xi)$ satisfies for any $\alpha \in [\alpha_2,1]$ and any  $\xi \in \Delta_{\mu_2}$
$$
\| \Rr_\alpha(\xi) \|_{\Ww^{k,1}_q(m^{-1})} \le 
\frac{C_3 + C_4 \, \| \Rr_1(\xi) \|_{W^{k+1,1}_{q+1}(m^{-1})}}
{1- C_5 \, (1-\alpha) \, \| \Rr_1(\xi) \|_{W^{k+1,1}_{q+1}(m^{-1})}} \, , 
$$
for some constants $C_i$, $i=3, 4, 5$.
\end{itemize}

\end{lem}

The rest of the proof is done the same, and we recall the main steps. 

\medskip\noindent
{\bf Step 1: Structure of the spectrum and rough localization}

\smallskip
We decompose
$$
\hat\LL_\alpha = T + K
$$
with $T= \tau_\alpha \, \Delta - L(\bar F_\alpha)$, 
$K = 2 \, \tilde Q^+_\alpha(\bar F_\alpha,\cdot) - \bar F_\alpha \, L (\cdot)$. 
Since the spectrum of $T$ is easily seen to be included in 
$\Delta_{\mu_2 + \OO(1-\alpha)}^c$ (see the definition of $\mu_2 < 0$ in the 
statement of Theorem~\ref{theoL1} and using the regularity estimates on the 
stationary solutions of Proposition~\ref{estimatesonGe}), and since 
the operator $K$ is $T$-compact (it is the only point where we use the sharp 
estimate ``in the norm of the graph'' stated in Propositions~\ref{ContQe2} \& \ref{LaL1}), 
we conclude (for a small enough inelasticity) thanks to Weyl's theorem 
that $\Sigma (\hat\LL_\alpha) \cap \Delta_{\bar\mu}$ 
only contains discrete spectrum (that is isolated and finite multiplicity spectrum values, 
or, in other words, eigenvalues), a set that we denote by $\Sigma_d (\hat\LL_\alpha)$. 
Moreover, thanks to point (iii) in Lemma \ref{techLaBdelta} and to the 
geometry of the spectrum of $\hat \LL_1$ as stated in Theorem~\ref{theoL1} 
(more particularly point (iii)) we deduce that $\Sigma (\hat\LL_\alpha) \cap \Delta_{\bar\mu}$ 
is confined to a disc $B(0,c \, (1-\alpha))$. 
To sum up, we have yet proved 
$$
\Sigma(\hat\LL_\alpha) \cap \Delta_{\bar\mu} \subset \Sigma_d(\hat\LL_\alpha)  \cap  B(0,c \, (1-\alpha)).
$$

\medskip\noindent
{\bf Step 2: The ``energy eigenvalue''}
\smallskip

First, we infer (for instance, from the decomposition (\ref{decomp1})) that any 
eigenfunction associated to an eigenvalue in $\Delta_{\mu_2}$ is smooth: 
more precisely, there exists $C = C_{k,q,m}$ such that $\forall \, (\lambda,\psi)$ 
with $\lambda \in \Delta_{\mu_2}$, $\psi \in \Ww^2_1$ satisfying 
$\hat \LL_\alpha \, \psi = \lambda \, \psi$, 
there holds $\| \psi \|_{W^{k,1}_q(m^{-1})} \le C \, \| \psi \|_{L^1_2}$. 

\smallskip
Second, we introduce the projector $\Pi_\alpha$ on the eigenspace associated 
to the spectrum included in $\Delta_{\mu_2}$. It is given by 
$$
\Pi_\alpha := - {1 \over 2 \, \pi \, i} \int_{\{ \zeta \in \CC, \, |\zeta|=r \} } \Rr_\alpha(\zeta) \, d\zeta 
$$
In particular the operator  $\Pi_1$ is the projection on the energy eigenline $\R \, \phi_1$, where $\phi_1$ is the energy eigenfunction defined in Theorem~\ref{theoL1}. Thanks to the smoothness estimate on eigenfunctions just mentioned above and the estimate on $\LL_1 - \LL_\alpha$ in Lemma~\ref{techLaBdelta} (i) we easily deduce that 
$\|Ê\Pi_1 - \Pi_\alpha \|Ê< 1$ for $\alpha$ close enough to $1$, and then that $\hbox{rang} \, \Pi_\alpha = \hbox{rang} \, \Pi_1$. We may sum up as
$$
\Sigma(\hat\LL_\alpha) \cap \Delta_{\bar\mu} =  \{ \mu_\alpha \}, \quad \mu_\alpha \in \R, \,\,  |\mu_\alpha |Ê\le c \, (1-\alpha),
\quad \hbox{Nul}(\mu_\alpha - \hat \LL_\alpha) = \R \, \phi_\alpha.
$$

\smallskip
Third, we prove that $\| \phi_\alpha - \phi_1 \|_{W^{k,1}_q(m^{-1})} \le C \, (1-\alpha)$ 
exactly as in~\cite{MMGranular3}.

\smallskip
Fourth and last, we obtain a first order expansion of $\mu_\alpha$ thanks to the following computation. By integrating the eigenvalue equation related $\mu_\alpha$ 
$$
\hat\LL_\alpha \phi_\alpha = \mu_\alpha \, \phi_\alpha
$$
against $|v|^2$ and dividing it by $(1-\alpha)$, we get
$$
\frac{\mu_{\alpha}}{1-\alpha}  \, \EE(\phi _{\alpha})= \rho \,  2 \, N \, \rho(\phi_{\alpha}) - 2 \, (1+\alpha) \, 
\tilde D ( \bar F_\alpha, \phi_{\alpha}).
$$
Using convergence of $\bar F_\alpha \to \bar F_1$ and 
$\phi_\alpha \to \phi_1$ established before, and the fact 
that the mass of $\phi_1$ is zero, we deduce that 
\beqn\label{vp2:1}
\frac{\mu_{\alpha}}{1-\alpha}  \, \EE(\phi _{1}) = 
 - 4 \, \tilde D ( \bar F_1, \phi_1) + \OO(1-\alpha).
\eeqn
Then we compute thanks to (\ref{Mv2}) and (\ref{Mv4})
\beqn\label{vp2:2}
\EE(\phi_1) = 2 \, N \, c_0 \, \rho \, \bar\theta_1^2, 
\eeqn
where $c_0$ is still the normalizing constant in (\ref{phiaTOphi1}) such that $\| \phi_1 \|_{L^1_2} = 1$. Similarly, using (\ref{MMu3}),  (\ref{MMv2u3}) and the relation (\ref{tempSS}) which make a link between $b_1$ and $\bar\theta_1$, 
we find
\beqn\label{vp2:3}
\tilde D ( \bar F_1, \phi_1) = {3 \over 2} \, N \, c_0 \, \rho^2 \, \bar \theta_1.
\eeqn
We conclude gathering (\ref{vp2:1}), (\ref{vp2:2}) and (\ref{vp2:3}).


\section{Convergence to the stationary solution} \label{sec:cvgce}
\setcounter{equation}{0}
\setcounter{theo}{0}

In this section, we consider the nonlinear evolution equation 
(\ref{eqBolmua}) and we prove the convergence of its solutions to 
the stationary solution. 

\subsection{The results}

We first state a local linearized stability result.

\begin{prop}\label{domaineunif} 
For any $\alpha \in [\alpha_3,1)$, the stationary solution $\bar F_\alpha$ 
is locally asymptotically stable, with domain of stability uniform according to $\alpha \in [\alpha_3,1)$. 

More precisely, let us fix $\rho \in (0,\infty)$ and some exponential weight function $m$ as in~(\ref{defdem}). 
There is $k_1,q_1 \in  \N^*$ such that for any $M_0 \in (0,\infty)$ there exists $C,\eps \in (0,\infty)$ such 
that for any $\alpha \in [\alpha_3,1]$, for any $f_{\mbox{\scriptsize{{\em in}}}} \in H^{k_1} \cap L^1(m^{-q_1})$ 
with mass $\rho$, momentum $0$ satisfying
\beqn\label{hypginGa}
\| f_{\mbox{\scriptsize{{\em in}}}} \|_{H^{k_1} \cap L^1(m^{-q_1})} \le M_0,
\qquad  \|f_{\mbox{\scriptsize{{\em in}}}} - \bar F_\alpha \|_{L^1(m^{-1})} \le \eps,  
\eeqn
the solution $f$ to the equation (\ref{eqBolmua}) with initial datum $f_{\mbox{\scriptsize{{\em in}}}}$ satisfies 
\beqn \label{energyGa}
\forall \,  t \ge 0, \quad \| \Pi_\alpha \, (f_t - \bar F_\alpha ) \|_{L^1(m^{-1})} 
\le C \, \| f_{\mbox{\scriptsize{{\em in}}}} - \bar F_\alpha \|_{L^1(m^{-1})} \, 
e^{\mu_\alpha  \, t},
\eeqn
\beqn\label{lingtoGa}
\forall \,  t \ge 0, \quad \| (\mbox{{\em Id}} - \Pi_\alpha) \, (f_t - \bar F_\alpha ) \|_{L^1(m^{-1})} 
\le C \, \| f_{\mbox{\scriptsize{{\em in}}}} - \bar F_\alpha \|_{L^1(m^{-1})} \, 
e^{(3/2) \, \mu_\alpha  \, t}.
\eeqn
\end{prop}

Then we prove that when the inelasticity is small, depending on the size of 
the initial datum (but not necessarily close to the stationary solution), the equation~(\ref{eqBolmua}) 
is stable. This mainly relies on the fact that the entropy production timescale 
is of a different order (much faster) that the energy dissipation timescale as $\alpha \to 1$.

  \begin{prop}\label{Attractiveg} 
  Define $k_2 := \max\{ k_0, k_1 \}$, $q_2 := \max \{ q_0,q_1, 3 \}$, where $k_i$ and $q_i$ are 
  defined in Theorem~\ref{CSK&EEP} and Corollary~\ref{DHeDH1}. For any $\rho, \, \EE_0, \, M_0$ 
  there exists $\alpha_4 \in [\alpha_3,1)$, $c_1 \in (0,\infty)$ and for any $\alpha \in [\alpha_4,1]$ 
  there exist $\varphi=\varphi(\alpha)$ with $\varphi(\alpha) \to 0$ as $\alpha \to 1$ 
  and $T=T(\alpha)$ (possibly blowing-up as $\alpha \to 1$) such that for any initial datum 
  $0 \le f_{\mbox{{\em \scriptsize{in}}}} \in L^1 _{q_2} \cap H^{k_2} \cap \Cc_{\rho,0,\EE_0}$ with 
  $$
  \| f_{\mbox{{\em \scriptsize{in}}}}\|_{ L^1 _{q_2} \cap H^{k_2} } \le M_0,
  $$
the solution $f$ associated to the rescaled equation~(\ref{eqresca}) satisfies
  $$
  \forall \, t \ge 0, \quad \EE (f_t) \ge c_1 
  $$
  and for all $\alpha' \in [\alpha_4,1)$ and then all $\alpha \in [\alpha',1]$
  \beqn\label{estimateonGe3}
  \forall \, t \ge T(\alpha'), \quad  \left\| f_t - \bar F_\alpha \right\|_{L^1 _2} \le \varphi(\alpha').
  \eeqn
\end{prop}

Then the proof of the gobal convergence for smooth initial data 
only amounts to connect the two previous results of Propositions~\ref{domaineunif} 
and~\ref{Attractiveg} by choosing $\alpha$ such that $\varphi(\alpha) \le \eps$ 
where $\eps$ is the size of the attraction domain in Proposition~\ref{domaineunif} and 
$\varphi(\alpha)$ is defined in  Propositions~\ref{Attractiveg}. More precisely, 
we have straightforwardly the

\begin{cor}\label{coro:cvgHk}
Let us fix an exponential weight function $m$ as in (\ref{defdem}), 
with exponent $s \in (0,1)$.  
Then for any $\rho, \, \EE_0, \, M_0$  there exists $C$ and $\alpha_5 \in [\alpha_4,1)$ 
(depending on $\rho, \, \EE_0, \, M_0,m$) 
such that for any $\alpha \in [\alpha_5,1)$ and any initial datum 
$0 \le f_{\mbox{{\em \scriptsize{in}}}} \in L^1(m^{-q_2}) \cap H^{k_2} $ satisfying 
$$
f_{\mbox{{\em \scriptsize{in}}}}\in \Cc_{\rho,0,\EE_0}, \qquad 
\| f_{\mbox{{\em \scriptsize{in}}}}\|_{ L^1(m^{-q_2}) \cap H^{k_2} } \le M_0,
$$
the solution $f$ associated to the rescaled equation~(\ref{eqresca}) satisfies 
$$ 
\forall \,  t \ge 0, \quad \| \Pi_\alpha \, (f_t - \bar F_\alpha ) \|_{L^1(m^{-1})} 
\le C \, e^{\mu_\alpha  \, t},
$$
$$ 
\forall \,  t \ge 0, \quad \| (\mbox{{\em Id}} - \Pi_\alpha) \, (f_t - \bar F_\alpha )  \|_{L^1(m^{-1})} 
\le C \, e^{(3/2) \, \mu_\alpha  \, t}.
$$
\end{cor}

As a by-product of the previous propositions, we state and prove a result which provides a 
partial answer to the question (important from the physical 
viewpoint) of finding Liapunov functionals for this particles system. 
Let us define the required objects. We consider a fixed mass $\rho$ and 
some restitution coefficient $\alpha$ whose range will be specified below. 
At initial times, non-linear effects dominate and therefore we define 
$$
{\mathcal H}_1 (f) := H(g|M[f]) + \big(\EE- \bar\EE_\alpha \big)^2
$$
where $\bar \EE_\alpha = \EE(\bar F_\alpha)$ is the energy of the self-similar profile corresponding to $\alpha$ 
and the mass $\rho$.
At eventual times, linearized effects dominate. Therefore we define a quite natural candidate 
from the spectral study: 
$$
{\mathcal H}_2(f) := \| h^1 \|_{L^1(m^{-1})}^2 
+ (1-\alpha) \int_0 ^{+\infty} \left\| \mathcal{R}_\alpha(s) \, h^2 \right\|^2 _{L^2} \, ds,
$$
with $h^1 = \Pi_\alpha h$, $h^2 = \Pi_\alpha^\perp h$ and $h = g - \bar G_\alpha$.

\begin{prop} \label{prop:lyap} 
There is $k_4 \in \N$ big enough (this value is specified in the proof) such that for any exponential weight 
function $m$ as defined in (\ref{defdem}), any time $t_0 \in (0,\infty)$ and any $\rho,\EE_0, M_0 \in (0,\infty)$, 
there exists $\kappa_* \in (0,\infty)$ and $\alpha_6 \in [\alpha_5,1)$ such that for any $\alpha \in [\alpha_6,1]$ 
and initial datum $f_{\mbox{\scriptsize{{\em in}}}} \in H^{k_4} \cap L^1(m^{-1})$ satisfying
$$ 
f_{\mbox{\scriptsize{{\em in}}}} \in \Cc_{\rho,0,\EE_0}, 
\qquad \| f_{\mbox{\scriptsize{{\em in}}}} \|_{H^{k_4} \cap L^1(m^{-1})} \le M_0, \qquad
f_{\mbox{\scriptsize{{\em in}}}}(v) \ge M_0 ^{-1} \, e^{-M_0 \, |v|^8}, 
$$
the solution $g$ to the rescaled equation (\ref{eqresca}) with initial datum 
$f_{\mbox{\scriptsize{{\em in}}}}$ is such that the functional 
$$
{\mathcal H}(f_t) = {\mathcal H}_1(f_t) \, {\bf 1}_{\big\{{\mathcal H}_1(f_t) \ge \kappa_*\big\}} 
+ {\mathcal H}_2(f_t) \, {\bf 1}_{\big\{{\mathcal H}_1(f_t) \le \kappa_*\big\}} 
$$
is decreasing for all times $t \in [0,+\infty)$. Moreover, $\HH(f(t, \cdot))$ is strictly decreasing as 
long as $f(t, \cdot)$ has not reached the self-similar state $\bar F_\alpha$. 
\end{prop}

Then the previous results can be extended to general initial data by a study of the decay 
of singularities in the same spirit as in~\cite{MMGranular3}. Indeed one can easily prove 
the following decomposition result.

 \begin{lem} \label{decompos2} 
 Consider $f_{\mbox{\scriptsize{{\em in}}}} \in L^1_3$ and the 
 associated solution $f \in C([0,\infty);L^1_3)$ to the rescaled equation (\ref{eqresca}). 
 Assume that for some constant $\rho, c_1,M_1,T \in (0,\infty)$ there holds
 \beqn\label{assumgdecomp1}
 f_{\mbox{\scriptsize{{\em in}}}} \in \Cc_{\rho,0}, \qquad 
 \| f_{\mbox{\scriptsize{{\em in}}}} \|_{L^1_3} \le M_1, \qquad \forall \, t \in [0,T], \quad 
 \EE(f(t,\cdot)) \ge c_1.
 \eeqn
Then, there are $\alpha_7 \in [\alpha_6,1)$ and  $\lambda \in (-\infty,0)$, 
and for any exponential weight function $m$ (as defined in (\ref{defdem}) 
and any $k \in \N$, there exists a constant $K$ (which depends on $\rho, c_1,M_1,k,m$) 
such that for any $\alpha \in [\alpha_7,1]$, we may split $f = f^S + f^R$ with 
\beqn\label{concludegdecomp1bis}
\forall \, t \in [0,T], \quad  \| f^S(t, \cdot) \|_{H^k \cap L^1(m^{-1})} \le K, \qquad 
\| f^R(t, \cdot) \|_{ L^1_3} \le K \, e^{\lambda \, t}.
\eeqn
 \end{lem}

This last tool together with the $L^1$ usual stability result of the elastic Boltzmann 
(stating that the error between two flows grows at most exponentially) allow to 
conclude to the proof of point (iv) of Theorem~\ref{theo:uniq}.

\subsection{Sketch of the proof}
\label{subsec:sketchnonlinear}

The proof is exactly the same as in~\cite[Section~6]{MMGranular3}. Indeed the precise form of the 
equation was not used in this proof, it only uses: 

\begin{itemize}
\item on the one hand, for the perturbative argument, the spectral properties obtained before (and valid 
here), some uniform estimates on moments and regularity in terms of a lower bound on the energy (still valid 
here) and some estimates from above on the bilinear term (which are the same here);  
\item on the other hand, for the argument ``in the large'', it uses elastic entropy - entropy production estimates 
(independent of our problem), and the difference of timescales between entropy production (which 
is $\OO(1)$) and temperature thermalization (which is $\OO(1-\alpha)$), also present here.
\end{itemize}


\appendix

\section{Appendix: Functional toolbox on the collision operator}
\setcounter{equation}{0}
\setcounter{theo}{0}

In this first appendix we shall recall some functional results on the collision operator $Q$ 
which were obtained in \cite{MMGranular3}.

First let us show that the collision operator depends continuously 
on the inelasticity coefficient $\alpha \in [0,1]$. 
Since it is an unbounded  operator, this continuous dependency is expressed 
in the norm of the graph of 
the operator or in some weaker norm. We start showing that this dependency of 
the collision operator is Lipschitz, and even $C^{1,\eta}$ for any $\eta \in (0,1)$, 
when allowing a loss (in terms of derivatives
and  weight) in the norm it is expressed. Let us define the formal derivative of the collision 
 operator according to $\alpha$ by 
 $$ 
 Q'_\alpha(g,f) :=  \nabla_v \cdot \left( 
 \int_{\R^N} \!\! \int_{\Sph^{N-1}} g('v_*(\alpha)) \, f('v(\alpha))
 \, b \,  |u| \, \left( {u-|u| \, \sigma  \over 4 \,  \alpha^2} \right) \, d\sigma \, dv_*\right)  
 $$
 or by duality 
 $$ 
 \langle Q'_\alpha(g,f) , \psi \rangle := 
 \int_{\R^N} \!\!  \int_{\R^N} \!\! \int_{\Sph^{N-1}} g_* \, f \, 
 b \,  |u| \, \left( {|u| \, \sigma - u  \over 4} \right)  \, 
 \nabla \psi (v'_\alpha) \,  d\sigma \, dv_* \, dv.
 $$

\begin{prop}[See Proposition~3.1 in \cite{MMGranular3}]\label{ContQe1} 
Let us fix a smooth exponential weight $m = \exp ( - a \, |v|^s )$, $a \in (0,+\infty)$, $s \in (0,1)$.
Then 
\begin{itemize}
\item[(i)] For any $k,q \in \N$ the exists $C \in (0,\infty)$ such that for any smooth functions $f, g$
(say in $\Ss(\R^N)$) and any $\alpha \in [0,1]$ there holds
\bear \label{QaWk1}
&&\big\|  Q^\pm_{\alpha}(g,f) \big\|_{W^{k,1}_q(m^{-1})} 
 \le C_{k,m} \, \|f \|_{W^{k,1} _{q+1} (m^{-1})} \, \| g \|_{W^{k,1} _{q+1} (m^{-1})} \\ \label{Q'aWk1}
 &&\big\|  Q'_{\alpha}(g,f) \big\|_{W^{k,1}_q(m^{-1})} 
 \le C_{k,m} \, \|f \|_{W^{k+1,1} _{q+2} (m^{-1})} \, \| g \|_{W^{k+1,1} _{q+2} (m^{-1})}.
 \eear
\item[(ii)] Moreover,  for any smooth functions $f, g$ and for any $\alpha, \alpha' \in [0,1]$, there holds
  \bear \nonumber
&& \quad \big\| Q^+_\alpha(g,f) - Q^+_{\alpha'}(g,f)- (\alpha-\alpha') \, Q'_\alpha(g,f) \big\|_{W^{-2,1}_q(m^{-1})}  \\ \label{Q'}
&& \qquad \qquad   \le |\alpha-\alpha'|^2\,    \|f \|_{L^1_{q+3} (m^{-1})} \,  \| g \|_{L^1_{q+3}(m^{-1})}.
 \eear
\item[(iii)] As a consequence, there holds
 \begin{equation}   \label{convLetoL1} 
 \mbox{ } \quad \big\| Q^+_{\alpha'}(g,f) - Q^+_{\alpha}(g,f) \big\|_{W^k _q(m^{-1})}
 \le C \, |\alpha-\alpha'| \,  \|f \|_{W^{2k +3,1}_{q+3} (m^{-1})} \,  
  \| g \|_{W^{2k +3,1}_{q+3} (m^{-1})},
 \end{equation}
and for any $\eta \in (1,2)$, there exists $k_\eta \in \N$, $q_\eta \in \N$ and $C_\eta \in (0,\infty)$ such that 
  \bear \nonumber
&& \quad \big\| Q^+_\alpha(g,f) - Q^+_{\alpha'}(g,f)- (\alpha-\alpha') \, Q'_\alpha(g,f) \big\|_{L^1(m^{-1})}  \\ \label{expLatoL1}
&& \qquad \qquad  
  \le C_\eta \, |\alpha-\alpha'|^\eta \, 
  \|f \|_{W^{k_\eta,1} _{q_\eta}(m^{-1})} \,  \| g \|_{W^{k_\eta,1}_{q_\eta}(m^{-1})}.
 \eear
\end{itemize}
\end{prop}

\medskip
We next state a mere (H\"older) continuity dependency on $\alpha$, which is however stronger than Proposition~\ref{ContQe1} in some sense, since it is  written in the norm of the graph of the operator for one of the arguments. 
 
\begin{prop}[See Proposition~3.2 in \cite{MMGranular3}]\label{ContQe2} 
For any $\alpha,\alpha' \in (0,1]$,  and any $g \in L^1_1(m^{-1})$, 
$f \in W^{1,1}_1 (m^{-1})$, there holds
 \beqn\label{convLetoL2}
 \quad \left\{
  \begin{array}{l}
  \big\|Q^+_\alpha(g,f) - Q^+_{\alpha'}(g,f) \big\|_{L^1(m^{-1})}
                 \le \eps (\alpha-\alpha') \, \| f \|_{W^{1,1}_{1}(m^{-1})} \, 
                 \|g \|_{L^1_1(m^{-1}) }, \vspace{0.3cm} \\
  \big\|Q^+_\alpha(f,g) - Q^+_{\alpha'} (f,g) \big\|_{L^1(m^{-1})} 
                 \le \eps (\alpha-\alpha') \, \| f \|_{W^{1,1}_{1}(m^{-1})} \, \|g  \|_{L^1_1(m^{-1}) }.
  \end{array}
  \right.
  \eeqn
where $\eps(r) = C \, r^{{1 \over 3 + 4 / s}}$ for some constant $C$ (depending only on $b$). 
\end{prop}

As a consequence of Propositions \ref{ContQe1} \& \ref{ContQe2}, 
together with Lemma~\ref{Contderivee=1Ge} and Proposition~\ref{estimatesonGe}, we have the

\begin{prop}[See Proposition~5.8 in \cite{MMGranular3}]\label{LaL1} 
For any $k,q \ge 0$ and for any exponential weight function $m$ there exists some constant $C$ such that $\forall \, \alpha \in (\alpha_0,1]$
\bean
&(i)&  \left\|  \LL^+ _\alpha \right\|_{W^{k,1}_q(m^{-1}) \to W^{k,1}_{q+1} (m^{-1})} \le C, \\
&(ii)&  \left\| \LL^+ _1 - \LL^+ _\alpha \right\|_{W^{k,1}_q(m^{-1}) \to W^{2\,k+3,1}_{q+3} (m^{-1})}  \le C \, (1-\alpha)\\
&(iii)&
 \left\| \LL^+ _1 - \LL^+ _\alpha \right\|_{W^{k,1}_q(m^{-1}) \to W^{k,1}_{q+1} (m^{-1})}  \le 
       \eps(1-\alpha),
\eean
where $\eps$ is defined as previously (up to a constant). 
\end{prop}

Second we shall give some estimates on the entropy production functional associated with the 
elastic collision operator $Q_1$. We begin with a simple consequence of Proposition~\ref{ContQe1}. 

\begin{cor}[See Corollary~3.4 in \cite{MMGranular3}]\label{DHeDH1} There exists $k_0, q_0 \in \N$ such that for any $a_i \in (0,\infty)$ $i=1, \, 2, \, 3$, 
there exists an explicit constant $C \in (0,\infty)$ such that for any function $g$ satisfying 
$$ 
 \| g \|_{H^{k_0} \cap L^1_{q_0}} \le a_1, \quad g \ge a_2 \, e^{-a_3 \, |v|^8},
$$
there holds
$$
\big| D_{H,\alpha}(g) - D_{H,1}(g) \big| \le C \, (1-\alpha),
$$
where we recall that $D_{H,\alpha}$ is defined by
\beqn\label{defDHa}
D_{H,\alpha} (g) = {1 \over 2} \int\!\! \int\!\! \int_{\R^{2N}\times S^{N-1}} g \, g_* \left( {g'g'_* \over g g_*} -\log {g'g'_* \over g g_*} 
     -  1 \right) \, B \, dv \, dv_* \, d\sigma \ge 0.
\eeqn
\end{cor}

\medskip
Let us recall now two famous inequalities, namely the Csisz\'ar-Kullback-Pinsker inequality (see \cite{Csisz,Kul59}) and
the so-called entropy-entropy production inequalities (the version we present here is established in \cite{VillCerc}) that we will use several time in the sequel. 

\smallskip
  \begin{theo}[See Theorem~3.5 in \cite{MMGranular3}]\label{CSK&EEP} 
\begin{itemize}
\item[(i)] For a given function $g \in L^1_2$, let us denote by $M[g]$ the Maxwellian function 
with the same mass, momentum and temperature as $g$.
For any $0 \le g \in L^1_2(\R^N)$, there holds
\beqn\label{CKineg}
\big\| g - M[g] \big\|_{L^1}^2 \le 2 \, \rho(g) \,  \int_{\R^N} g \ln { g \over M[g]} \, dv.
\eeqn
\item[(ii)] For any $\eps > 0$ there exists $k_\eps, \, q_\eps \in \N$  and for any $A \in (0,\infty)$ there exists $C_\eps = C_{\eps,A} \in (0,\infty)$ such that for any $g \in H^{k_\eps} \cap L^1_{q_\eps}$ such that 
$$ 
g(v) \ge A^{-1} \, e^{-A \, |v|^8}, \quad \| g \|_{ H^{k_\eps} \cap L^1_{q_\eps}} \le A,
$$
there holds
\beqn\label{EEPineg}
C_\eps \, \rho(g)^{1-\eps} \, \left( \int_{\R^N} g \ln { g \over M[g]} \, dv \right)^{1+\eps} \le D_{H,1}(g).
\eeqn
\end{itemize}
\end{theo}

\medskip
Third and last we recall some results on the linearized elastic collision operator.

\begin{theo}[See \cite{GM:04} and Proposition~5.5 \& 5.7 in \cite{MMGranular3}]\label{theoL1} 

Let define $\hat\LL_1$ as the restriction of $\LL_1$ to the space $L^1_1 \cap \Cc_{0,0}$. Then for any $k,q \ge 0$ and for any exponential weight function $m$ the following holds.

\begin{itemize}
\item[(i)]  The regularized truncation $\LL^+ _{1,\delta}$ of the ``gain'' part (see (\ref{decomp2}) and \cite{GM:04}) satisfies
\[ 
\left\| \LL^+ _1 - \LL^+ _{1,\delta} \right\|_{W^{k,1}_q(m^{-1}) \to W^{k,1} _{q+1} (m^{-1}) } \le \eps (\delta) 
\]
 where $\eps (\delta) >0$ is an explicit constant going to $0$ as $\delta$ goes to $0$. For any $\delta >0$, the linear operator  $A_\delta$ (see (\ref{decomp1})) satisfies $A_\delta: L^1 \to W^{\infty,1} _\infty(m^{-1})$ is bounded.
 
\item[(ii)] The unbounded operator  $\hat\LL_1$ is well defined and close in $\Ww^{k,1} _q(m^{-1})$. Its spectrum $\Sigma(\hat  \LL_1)$ is real and satisfies
$$
\Sigma(\hat \LL_1) \subset \Delta_{\mu_2} ^c \cup \{ \mu_1 \},
$$
with $\mu_2 < 0$ and $\mu_1=0$  is the $1$-dimensional "energy eigenvalue" associated to the  "energy eigenfunction" $\phi_1 := c_0 \, (|v|^2 - N \, \bar\theta_1) \, \bar F_1$. In particular,  $\hat \LL_1$ is onto from $\OO \cap \Ll^1_1(m^{-1})$ onto $\OO$.

\item[(iii)] The resolvent   $\Rr_1(\xi)$ has a sectorial property for the spectrum substracted from the 
``energy'' eigenvalue, namely there is a constructive $\lambda \in (\mu_2,0)$ such that  
$$
\forall \, \xi \in \Aa, \quad \| \Rr_1 (\xi) \|_{\Ww^{k,1}_q(m^{-1})} \le a_{k,q} + \frac{b_{k,q}}{|\xi + \lambda|}, 
$$
with some explicit constant $a_{k,q}, b_{k,q} >0$ and
$$
\Aa = \left\{ \xi \in \mathbb{C}, \quad \mbox{{\em arg}}(\xi +\lambda) 
\in \left[-\frac{3 \pi}{4}, \frac{3 \pi}{4} \right] \ \mbox{ and } \ \Re e \, \xi \le \frac{\lambda}2 \right\}.
$$
\end{itemize}
\end{theo}

A quite simple consequence of point (ii) above is the following quantitative invertibility result on a modified version of the linearized elastic collision operator.

\begin{lem}[See Lemma~4.3 in \cite{MMGranular3}]\label{lemDPhiM} The linear operator $A \ : \ \Ll^1_1(m^{-1}) \to  \R \times  \OO $ defined thanks to (\ref{D2Phi}) is invertible: it is bijective with $A^{-1}$ bounded with explicit estimate. 
\end{lem}


\section{Appendix: Moments of Gaussians}
\setcounter{equation}{0}
\setcounter{theo}{0}

We state here some results on the moments of tensor product of Gaussians (the proof is done 
is the appendix of~\cite{MMGranular3}).

\begin{lem}\label{MomentsGausians} The following identities hold
\bear
\label{Mv2}
&&\int_{\R^N} M_{1,0,1} \, |v|^2 \, dv = N,  \\
\label{Mv4}
&& \int_{\R^N} M_{1,0,1} \, |v|^4 \, dv = N \, (N+2), \\
\label{MMu3}
&&\int_{\R^N \times \R^N} M_{1,0,1} \, (M_{1,0,1})_* \, |u|^3 \, dv \, dv_* = 
     2^{3/2} \, \int_{\R^N} M_{1,0,1} \, |v|^3 \, dv, \\ \label{MMv2u3}
&&\int_{\R^N \times \R^N} M_{1,0,1} \, (M_{1,0,1})_* \,|v|^2 \,  |u|^3 \, dv \, dv_* = 
 \sqrt{2} \, (2N+3) \, \int_{\R^N} M_{1,0,1}(v) \, |v|^3 \, dv.
\eear    
\end{lem}


\section{Appendix: Interpolation inequalities}
\setcounter{equation}{0}
\setcounter{theo}{0}

Again the proof of the following simple interpolation inequality can be found 
in the appendix of \cite{MMGranular3}.

\begin{lem} \label{interpolineg} 
\begin{itemize}
\item[(i)] 
For any $k,k^*,q,q^* \in \ZZ$ with $k\ge k^*$, $q\ge q^*$ and any $\theta \in (0,1)$ there is  $C \in (0,\infty)$ such that for 
$h \in W^{k^{**},1}_{q^{**}}(m^{-1})$
\beqn\label{hWk1qm}
\| h \|_{W^{k,1}_q(m^{-1})} \le C \,\| h \|_{W^{k^*,1}_{q^*}(m^{-1})}^{1-\theta} \,  \| h \|_{W^{k^{**},1}_{q^{**}}(m^{-1})}^\theta.
\eeqn
with $k^{**},q^{**} \in \ZZ$ such that $k = (1-\theta) \, k^{*} + \theta \, k^{**}$, $q = (1-\theta) \, q^{*} + \theta \, q^{**}$.
\item[(ii)] 
For any $k, q \in N^*$ and any exponential weight function $m$ as defined in (\ref{defdem}), 
there exists $C \in (0,\infty)$ such that for any $h \in H^{k^\ddagger} \cap L^1(m^{-12})$ 
with $k^\ddagger := 8 k + 7 (1+N/2)$ 
\beqn\label{hWk1qm2}
\| h \|_{W^{k,1}_q(m^{-1})} \le C \, \| h \|_{H^{k^\ddagger}}^{1/4} \, \,\| h \|_{L^1( m^{-12})}^{1/4} \,  \| h \|_{L^1(m^{-1})}^{3/4}.
\eeqn
\end{itemize}
\end{lem}


\section{Appendix: Lower bound for the diffusive inelastic Boltzmann equation}
\label{app:lb}
\setcounter{equation}{0}
\setcounter{theo}{0}

In this last appendix, we state and prove a technical result on the pointwise 
lower bound of the diffusive Boltzmann equation (the new difficulty as compared 
to previous results is due to the diffusion and not transport nature of the term 
added to the collision operator).

\begin{prop}\label{LowerBdga} Let $g \in C([0,\infty);L^1_3)$ be a solution of the rescaled evolution equation
$$
\partial_t g = Q(g,g) + \tau_\alpha \, \Delta_v g,
$$
with inelasticity parameter $\alpha \in (0,1)$, and assume that for some $C,T \in (0,\infty)$
$$
\sup_{[0,T]} \| g \|_{L^2 \cap L^1_3} \le C.
$$
For any $t_1\in (0,T)$ there exists $a_1 \in (0,\infty)$ (depending on $C$, $\rho$ and $t_1$ but not on $T$) such that 
\beqn\label{concluLowerBdga}
 \forall \, t \in [t_1,T], \,\,\, \forall \, v \in \R^N, \quad g(t,v) \ge a_1^{-1} \, e^{-a_1 \, |v|^8}.
\eeqn
\end{prop}

 We closely follow the proof of the Maxwellian lower bound for the solutions of the elastic Boltzmann equation 
(see \cite{Ca32,PW97}) taking advantage of some technical results established in its extension to  the solutions of 
the inelastic Boltzmann equation (see \cite[Theorem 4.9]{MMII} and \cite[Lemma 2.6]{MMGranular3}).
The starting point is again the evolution equation satisfied by $g$  written in the form
\beqn\label{evolutionEq} \quad
\partial_t g  -  \tau_\alpha \, \Delta_v g + \lambda(v)  \, g = 
Q^+_\alpha(g,g) + \big(\lambda(v) - L(g)\big)  \, g, \quad \lambda(v) = \kappa \, (1 +  |v|) ,
\eeqn
where the last term in the right hand side term is clearly non-negative for some well-chosen 
numerical constant $\kappa \in (0,\infty)$. 
Let us introduce the semigroup $S_t$ associated to the linear evolution equation $\partial_t g - \tau_\alpha \, \Delta_v \, g + \lambda(v) \, g = 0$.

\medskip\noindent
We start establishing some technical results that we need in the proof of Proposition~\ref{LowerBdga}. Let us recall some elementary results extracted from \cite{MMII} (which are a mere adaptation to the inelastic collision of some result proved in \cite{PW97}).

\begin{lem}\cite[Lemma 4.6]{MMII}\label{lower1} Let $0 \le \ell \in (L^1_2 \cap L^2)(\R^N)$ satisfy 
$$
\int_{\R^N} \ell \, dv = m_0, \quad \int_{\R^N} \ell \, |v|^2 \, dv \le m_1, \quad 
 \int_{\R^N} \ell^2\, dv \le m_2,
$$
for some positive real constants $m_i$. There exists $R > r > 0$ and $\eta > 0$ depending only on $m_0, m_1, m_2$, and $(v_i)_{i=1,\dots,4}$ such that $|v_i|\le R$, $i=1,\dots,4$, $|v_i - v_j|\ge 3r$ for $1\le i \not= j \le 3$, and
\beqn \label{intBrg}
\int_{B(v_i,r)} \ell (v) \, dv \ge \eta \quad\hbox{ for }\quad i=1,\,2,\,3,
\eeqn
\beqn \label{sphereplan}
\forall \, w_i \in B(v_i,r), \quad E^\alpha_{w_3,w_4} \cap S^\alpha_{w_1,w_2}
\hbox{ is a sphere of radius larger than } r,
\eeqn
where, as in \cite[Proposition 1.5]{MMII}, $E^\alpha_{w_3,w_4}$ stands for the
hyperplan orthogonal to the vector $w_3-w_4$ and passing through the point $\Omega(w_3, w_4)$, defined by 
\[ 
\Omega(w_3, w_4) := w_3 + (1- \alpha) \, (w_3-w_4) /(1+\alpha),
\] 
and $S^\alpha_{w_1,w_2}$ stands for the sphere of all possibles post-collisional 
velocity $v'$ defined by (\ref{vprimvprim*}) from $(v=w_3,v_*=w_4,\sigma)$.
\end{lem}

\begin{lem}\label{lower2} 
Let $0 \le f,g,h \in L^1(\R^N)$ and  $(v_i)_{i=1,\dots,4}$ 
satisfy (\ref{intBrg}), (\ref{sphereplan}) and $|v_i|\le R$, $i=1,\dots,4$, $|v_i - v_j|\ge 3r$ 
for $1\le i \not= j \le 3$, for some given constants $R > r > 0$, $\eta > 0$. 
There exists $T'_0 > 0$, $\delta'_0 > 0$ and $\eta'_0 > 0$ only depending on $R$, $r$ and $\eta$ such that 
$$
\forall \, t' \in [0,T_0] \qquad Q^+(f, S_{t'} Q^+ (g,h)) \ge \eta_0' \, {\bf 1}_{B(v_3, \delta'_0)}.
$$
\end{lem}

\smallskip\noindent
{\sl Proof of Lemma~\ref{lower2}}.  We first recall a convenient  
formula to handle representations of the iterated 
gain term. Using Carleman representation \cite[Proposition 1.5]{MMII}, 
for any $f$, $h$ and $\ell$ and any $v \in \RR^N$ there holds (setting $'v = w$ and $'v_* = w_*$)
\bean
&&Q^+(f,S_t \, Q^+(h,\ell)) (v)
= C'_b \int_{\RR^N}  {f(w) \over |v-w|} \left\{ \int_{E _{v,w}}
S_t \, Q^+(g,h)(w_*)  \, dw_* \right\}\, dw \\
&& \qquad =   \int_{\RR^{3N}} f(w)  \,  g(z) \, h(z_*) \,  \left\{ C'_b \,  { |z_* - z| \over |v-w|} 
\int_{S^{N-1}} (S_t \, {\bf 1}_{E^\alpha  _{v,w}}) (z') \, b(\sigma \cdot \hat z) \, d\sigma \right\} \, dz_* \, dz \, dw,
\eean
where $E_{v,w} ^\alpha$ is defined in the statement of Lemma~\ref{lower1} and $z'$ is defined from (\ref{vprimvprim*}) 
with $(v=z,v_*=z_*,\sigma)$. 
Let us define $\bar f  := f \, {\bf 1}_{B(v_4,r)}$, 
$\bar g = g \, {\bf 1}_{B(v_1,r)}$, $\bar h = h \, {\bf 1}_{B(v_2,r)}$ and $\tilde S_t$ the semigroup 
associated to evolution equation $\partial_t g -\tau_\alpha \, \Delta_v g - \lambda(2R) g = 0$. 
By the maximum principle, we have for any $v \in B(v_3,r)$ 
\bean
&&Q^+(f,S_t \, Q^+(h,\ell)) (v)
\ge   \int_{\RR^{3N}} \bar f(w)  \,  \bar g(z) \, \bar h(z_*) \,  \left\{ C''_b \,  { r \over R} 
\int_{S^{N-1}} (\tilde S_t \, \chi) (z') \, d\sigma \right\} \, dz_* \, dz \, dw.
\eean
Taking $v \in B(v_3,r)$, $w \in B(v_4,r)$, $z \in B(v_1,r)$, $z_* \in B(v_2,r)$ and denoting by $A$ the term between brackets, we have thanks to (\ref{sphereplan}) 
$$
A(v,w,z,z_*) := C''_b \,  { r \over R} \int_{S^{N-1}} \int_{\R^N} \chi(z'-u) \, e^{- {|u|^2 \over 2 \, \tau_\alpha \, t}} \, d\sigma du \ge  {r \over R} \,
C_b \, C \, r^{N-2} \, {1 \over 2},
$$
for any $\alpha \in [0,1]$ and $t \in [0,T_0]$ with $T_0$ small enough. We then conclude the proof. \qed

\medskip
For a given $\mu \in (0,1)$ we consider a function $ \chi \in W^{2,\infty}(\R_+)$, which satisfies the following properties: $\chi \ge 0$, $\chi' \le 0$, $\chi \equiv 1$ on $[0,\mu]$ and $\chi(x) = (1-x)^2$  for any $x \in [\mu',1]$ with $\mu' := \max\{1-1/N,(1+\mu)/2\}$. Abusing notations we define the radial symmetric functions $\chi$ and $\chid$ on $\R^N$ by setting $\chi(v) = \chi(|v|)$ and $\chid(v) = \chi(|v|/\delta)$ for any $\delta > 0$. 
 
\begin{lem}\label{lower3} For any $\delta_0 > 0$, there exists $C \in (0,\infty)$ such that 
$$
\forall \, \alpha \in (0,1), \,\,\,  \forall \, R > 0, \,\,\,  \forall \, \bar v \in B(0,R),
\,\,\,  \forall \, \delta > \delta_0, \,\,\, \forall \, t \ge 0,  \quad 
S_t \, (\tau_{\bar v} \chi_{_\delta}) \ge e^{-C \, (1 + R + \delta) \, t} \, (\tau_{\bar v} \chid).
$$
\end{lem}

\medskip
\noindent{\sl Proof of Lemma \ref{lower3}.} 
We shall rely on a maximum principle argument for the operator
$$
T f := \tau_\alpha \,  \Delta_v f - \lambda(v) \, f.
$$ 

First it can be seen easily that if $f_0$ is a non-negative function 
and $f=f(t,v)$ satisfies the differential inequality
$$
\partial_t f \ge T f, \quad t \in [0,T),
$$
then $f_t$ remains non-negative on this time interval $[0,T)$. 

Second let us show that
$$
\phi(t,v) := e^{-a \, t} \, \tau_{\bar v} \chi_\delta
$$
is a sub-solution for $a=a_{\delta,R}$ big enough, on a small initial 
time interval $[0,T)$, in the sense that
$$
\partial_t \phi \le T \phi, \quad t \in [0,T).
$$
After elementary computations and taking advantage of the fact 
that $\chi$ is radially symmetric, it amounts to show that 
\beqn\label{eqsubsol}
\big[a - \lambda(\delta \, w+\bar v)\big] \, \chi  
+ {\tau_\alpha \over \delta^2} \, \left[ \chi''  + \frac{(N-1)}{|w|} \, \chi' \right] 
\ge 0 \quad\hbox{on}\quad \R^N
\eeqn
for $a$ big enough. From the definition of $\chi$, it is zero outside of $B(0,1)$ and therefore equation (\ref{eqsubsol}) holds on $\R^N \backslash B(0,1)$. Observing that 
$$
\chi''(x)  + \frac{(N-1)}{x} \, \chi' (x) \ge 0
$$ 
for any $x \in [\mu', 1]$, we deduce that  equation  (\ref{eqsubsol}) holds on $ B(0,1)  \backslash B(0,\mu')$ as soon as $a \ge \kappa \, (1 + \delta + R)$. Finally, using that $\chi$ is decreasing and bounded in $W^{2,\infty}(0,1)$, we see that equation (\ref{eqsubsol}) holds on $B(0,\mu')$ as soon as 
$$
\big[a - \kappa \, (1+\delta+R)\big] \, \chi(\mu') \ge {\rho \over \delta^2} \,  \left\| \chi''(x)  + \frac{(N-1)}{x} \, \chi' (x) \right\|_{L^\infty(0,1)},
$$
and therefore as soon as $a \ge C \, (1+\delta+R)$ for some constant $C \in (0,\infty)$ only depending on $\chi$, $\delta_0$, $\rho$ and $N$.   
\qed

\medskip
\noindent{\sl Proof of Proposition \ref{LowerBdga}.} 
We split the  proof into two steps.

\smallskip\noindent
{\sl Step 1.} Thanks to the Duhamel formula a solution $g$ to the evolution equation (\ref{evolutionEq}) satisfies
\beqn\label{DuhamelF}
g(t,\cdot) \ge S_t \, g(0,\cdot) + \int_0^t S_{t-s} Q^+(g(s,\cdot),g(s,\cdot)) \,ds.
\eeqn
Let us fix $t_0 > 0$ and define $\tilde g_0 (t,\cdot) := 
g(t_0 + t, \cdot)$. Using twice the inequality  (\ref{DuhamelF}), we find 
\bear\label{glower2}
\tilde g_0(t,\cdot) &\ge&  \int_0^t  \int_0^s S_{t-s} Q^+\left(S_s \tilde g_0,
S_{s-s'} Q^+(S_{s'} \tilde g_0, S_{s'} \tilde g_0)  \right) \, ds' \, ds.
\eear
We first apply Lemma~\ref{lower2} 
to $\ell = S_{s'} \, \tilde g_0$ with $\tau = s-s'$ and we 
next  use Lemma~\ref{lower3} to obtain for any $t \in [0,T_0]$, 
with $T_0$ given by Lemma~\ref{lower2} 
\bear
\tilde g_0(t,\cdot) 
&\ge& \int_0^t \! \int_0^s S_{t-s} \eta \, {\bf 1}_{B(\bar v,r)} \, ds' \, ds \ge  \eta \,  \int_0^t  e^{-a \, (1- (t-s))} \, \tau_{\bar v} \chi_{_r}  \, s \, ds. 
\eear
We have then proved that there exists $T_1 > 0$ and $\bar{v} \in B(0,R)$ and for any $t_1 \in  (0,T_1/2]$ there exists $\eta_1 > 0$ such that 
\[
\forall \, t \in [0,T_1/2], \quad 
\tilde g_1( t, \cdot) := \tilde g_0( t+ t_1, \cdot) \ge \eta_1 \,  \tau_{\bar v} \chi_{_{\delta_1}} =: \eta_1 \, \bar\chi_{\delta_1}.
\]

\smallskip\noindent
{\sl Step 2. }  Using again the  Duhamel formula (\ref{DuhamelF}) and 
the preceding step we have
\[ 
\tilde g_1 (t,\cdot) \ge  \int_0^t S_{t-s} Q^+(\tilde g_1 (s,\cdot),\tilde g_1 
(s,\cdot)) \,  ds.
\]
Now let recall, that on the one hand, from  \cite[ Lemma 4.8]{MMII}, there exists $\kappa \in (0,\infty)$ such that 
$$
Q^+_\alpha ({\bf 1}_{B(0,1)},{\bf 1}_{B(0,1)}) \ge \kappa' \, {\bf 1}_{B(0,\sqrt{5}/2)}
$$
and that, on the other hand, the scaling properties of the $Q^+_\alpha$ term infers that 
$$
Q^+_\alpha(\phi(./\delta),\phi(./\delta)) (v) = \delta^{-N-1} \, Q^+_\alpha(\phi,\phi) (v/\delta)
$$ for any function $\phi$ and scaling coefficient $\delta > 0$. 
Then, thanks to Lemma~\ref{lower3}, there holds
\bean 
\tilde g_1 (t,\cdot)
&\ge&  \eta_1^2 \int_0^t S_{t-s} Q^+\big(\bar \chi_{_{\delta_1}},\bar \chi_{_{\delta_1}} \big) \,  ds \\
&\ge&  \eta_1^2 \, b \, e^{- a \,  (1 + R + \delta_1) \, t} \, t \, \delta_1^{-N-1} \, \bar \chi_{_{\theta \, \delta_1}}
\eean
on $[0,T_2]$ with $T_2 \in  (0,T_1/2]$, $\theta \in (1,\sqrt{5}/2)$ as close as we wish to $\sqrt{5}/2$ by choosing $\mu$ close to $1$ and $b$ a numerical cosntant (depending on $\kappa'$, $\chi$).  Defining $t_k = t_{k-1} + t_1 \, 2^{-k}$ and repeating the precedent computation we see that 
$$
g(t,\cdot) \ge \eta_k \, \bar  \chi_{_{\delta_k}} \qquad \forall \, t \in [t_k, T],
$$
with $\delta_k = \theta^k \, \delta_1$ and 
$$
\eta_{k+1} =  \eta_k^2 \, b \, e^{- a \,  (1 + R + \theta^k \, \delta_1) \, t_1 \, 2^{-k}} \, {t_1 \over 2^k} \, \delta_1^{-N-1} \, \left(\theta^{-N-1}\right)^k =  \eta_k^2 \, A \, B^{k},
$$
for some constants $A,B \in (0,\infty)$. Elementary computations yields to 
$$
\eta_{k} \ge \eta_1^{2^k} \, A^{1+2+...+2^{k-1}} \, B^{k+(k-1)\,2 + ...+2^{k-1}} \ge D^{2^k},
$$
for some constant $D \in (0,\infty)$. 
As a conclusion, using that $\theta^{8} > 2$, we have proved
$$ 
\forall \, t \ge t_1, \ \forall \, k \in \N, \quad g(t,v) \ge D^{2^k} \, {\bf 1}_{B(\bar v,2^{k/8} \, (\delta_1/2))}(v),
$$
from which we easily conclude. \qed

\bigskip
\noindent
{\bf{Acknowledgments.}}  The authors would like to thank Cyril Imbert for fruitful discussions 
related to Appendix~\ref{app:lb}.

\smallskip


\end{document}